\documentclass{article}

\usepackage{arxiv}

\usepackage{amssymb} %%%
\usepackage{amsmath} %%%%
\usepackage{amsthm}

\usepackage{subcaption} %%%%
\usepackage[export]{adjustbox} %%%%%

\usepackage[utf8]{inputenc} % allow utf-8 input
\usepackage[T1]{fontenc}    % use 8-bit T1 fonts
\usepackage{hyperref}       % hyperlinks
\usepackage{url}            % simple URL typesetting
\usepackage{booktabs}       % professional-quality tables
\usepackage{amsfonts}       % blackboard math symbols
\usepackage{nicefrac}       % compact symbols for 1/2, etc.
\usepackage{microtype}      % microtypography
\usepackage{lipsum}		% Can be removed after putting your text content
\usepackage{graphicx}
\usepackage{natbib}
\usepackage{doi}

\newtheorem{theorem}{Theorem}%  meant for continuous numbers
\newtheorem{proposition}[theorem]{Proposition}% 
\newtheorem{lemma}{Lemma}
\newtheorem{remark}{Remark}%

\title{A posteriori error analysis and adaptivity for a VEM discretization of the Navier-Stokes equations}

\author{{Claudio Canuto} \\
	Department of Mathematical Sciences\\
	Politecnico di Torino\\
	{Corso Duca degli Abruzzi 24}, {Turin}, {10129}, {Italy} \\
	\texttt{claudio.canuto@polito.it} \\
	\And
	{Davide Rosso} \\
	Department of Mathematical Sciences\\
	Politecnico di Torino\\
	{Corso Duca degli Abruzzi 24}, {Turin}, {10129}, {Italy} \\
	\texttt{daviderosso997@gmail.com} \\
}
\date{}

\renewcommand{\headeright}{}
\renewcommand{\undertitle}{}
\renewcommand{\shorttitle}{Adaptive VEMs for the Navier-Stokes problem}

\hypersetup{
pdftitle={A posteriori error analysis and adaptivity for a VEM discretization of the Navier-Stokes equations},
pdfsubject={},
pdfauthor={Claudio Canuto, Davide Rosso},
pdfkeywords={A posteriori estimator, Adaptivity, Virtual element method, Navier-Stokes equations, Computational fluid dynamics},
}

\begin{document}
\maketitle
\begin{abstract}
We consider the Virtual Element method (VEM)  introduced by Beir\~ao da Veiga, Lovadina and Vacca in 2016 for the numerical solution of the steady, incompressible Navier-Stokes equations; the method has arbitrary order $k \geq 2$ and guarantees divergence-free velocities. For such discretization, we develop a residual-based a posteriori error estimator, which is a combination of standard terms in VEM analysis (residual terms, data oscillation, and VEM stabilization), plus some other terms originated by the VEM discretization of the nonlinear convective term. 
We show that a linear combination of the velocity and pressure errors is upper-bounded by a multiple of the estimator (reliability). We also establish some efficiency results, involving lower bounds of the error. Some numerical tests illustrate the performance of the estimator and of its components while refining the mesh uniformly, yielding the expected decay rate. At last, we apply an adaptive mesh refinement strategy to the computation of the low-Reynolds flow around a square cylinder inside a channel.
\end{abstract}

% keywords can be removed

\section{Introduction}\label{sec1}

Since their inception (\cite{Basic,Guide}), Virtual Element Methods (VEM) have experienced an explosive growth of interest in the community devoted to the numerical solution of partial differential equations. VEM can be viewed as a generalization of Finite Element Methods, in which domain partitions formed by elements with rather arbitrary shapes are allowed; polytopal meshes are the prime example. The success of VEM is certainly due to the fact that they permit a high geometric flexibility (which, e.g., facilitates mesh generation, or local adaptive mesh refinements) without giving up the benefits of FEM, such as the use of weak formulations, or the use of polynomials in the calculation of elemental matrices.

Among the various applications of Virtual Element Methods, we highlight here the numerical solution of fluid-dynamics problems, such as the equations describing an incompressible flow. We refer to the papers \cite{Stokes_1, Navier-Stokes_1} which initiated this research line, followed by other contributions such as, e.g., \cite{a_post_Stokes,WANGNS}. A remarkable advantage of VEM over FEM is the easiness in generating discretizations which produce divergence-free velocities, thereby satisfying the continuity equation exactly. Another interesting feature is the possibility of handling meshes with cut elements, as those encountered, e.g., in problems of fluid-structure interactions; see \cite{Leaflet} for an example. Furthermore, adaptive mesh refinement -- a natural procedure in the presence of localized structures such as vortices -- gains efficiency from using virtual elements, for which the hanging nodes generated by local refinement are simply viewed as edge separators, and need not be removed. We refer to \cite{new_stab} for an interesting realization of this new perspective.

The present paper is a contribution to the development of adaptive Virtual Element Methods for the stationary incompressible Navier-Stokes equations in dimension 2. More precisely, we devise a residual-based a posteriori error estimator for the VEM discretization of order $k \geq 2$ introduced in \cite{Navier-Stokes_1}, we prove reliability and efficiency results, and we use the estimator in designing an adaptive refinement algorithm, which we apply to a classical fluid-dynamic problem in Wind Engineering, the flow around a square cylinder (\cite{Fluid_1}). 

Concerning the a posteriori error analysis for VEM, after the pioneering work \cite{a_post_cangiani} for a general second-order elliptic problem, we mention the a posteriori estimator given in \cite{a_post_Stokes} for the Stokes equations, and the one in \cite{WANGNS} for the Navier-Stokes equations. Note however that the VEM discretization considered hereafter differs from the one in \cite{WANGNS}, and we base our robustness analysis of the estimator on a general framework due to \cite{Ains}.

Our adaptive algorithm is based on the classical loop $\texttt{SOLVE} \to \texttt{ESTIMATE} \to \texttt{MARK} \to \texttt{REFINE}$, with $\texttt{SOLVE}$ realized by the VEM discretization \cite{Navier-Stokes_1} of order $k = 2$, and $\texttt{MARK}$ realized by D\"orfler's criterion. In the application, following \cite{Fluid_1}, we consider a 2D rectangular channel with a square obstacle symmetrically placed astride the horizontal axis, and we use meshes made of geometrically square elements, obtained by quad-tree refinements. These elements are indeed virtual elements: a hanging nodes created by a refinement splits the edge of a square into two aligned edges, thereby generating a polygon with 5 or more edges. 

The paper is organized as follows. Section \ref{sec:cont} introduces the steady continuous Navier-Stokes problem, while Section \ref{sec:discr} describes the VEM discretization adopted following \cite{Navier-Stokes_1}. Our main results related to the a posteriori error analysis is contained in Section \ref{sec:main}, in particular here we prove an upper bound of the error in terms of the estimator, together with some lower bounds. In Section \ref{sec:experiments}, we first validate our estimator with some numerical tests; then, we introduce the adaptive algorithm and we apply it to the flow around a square obstacle at Reynolds number $< 60$ (for which the solution is steady), showing the sequence of adaptive meshes and producing an accurate estimate of the length of the recirculation zone past the obstacle. The conclusions follow in Section \ref{sec:conclusions}.

\smallskip
{\em Notation.} We denote by $\|\cdot\|_{0,A}$ the $L^2(A)$-norm of a function (which may be a scalar, a vector, or a tensor) defined on a Lipschitz domain $A$ (in one or two dimensions), and by $( \cdot \, , \cdot)_A$ the inner product in $L^2(A)$. Similarly, we denote by $\vert \cdot \vert_{k,A}= \vert \cdot \vert_{H^k(A)}$ and $\Vert \cdot \Vert_{k,A}= \Vert \cdot \Vert_{H^k(A)}$ the Sobolev semi-norm and norm of order $k$ in $A$.  Furthermore, the notation $P \lesssim Q$ means $P \leq C \, Q$, where $C>0$ is a constant independent of the considered discretization and the Reynolds number.

\section{The continuous problem}\label{sec:cont}
We consider the steady Navier-Stokes problem in a polygonal domain $\Omega\subset\mathbb{R}^2$ under mixed Dirichlet-Neumann boundary conditions:
\begin{equation}\label{NS_forte}
    \text{\em find} \quad (\mathbf{u},p) \quad \text{\em such that} \quad
    \begin{cases}
    -\nu\mathbf{\Delta}\mathbf{u}+(\mathbf{\nabla}\mathbf{u})\mathbf{u}-\nabla p =\mathbf{f} \quad \text{in $\Omega$}\,, \\
    \nabla \cdot \mathbf{u}=0 \quad \text{in}\quad \Omega\,, \\
    \mathbf{u}=\mathbf{g} \quad \text{on} \quad \Gamma_D \subseteq \partial\Omega\,, \\
    (\nu\mathbf{\nabla}\mathbf{u}+p\mathbf{I})\mathbf{n}=\mathbf{t} \quad \text{on} \quad \Gamma_N = \partial\Omega \setminus \Gamma_F \,,\\
    \end{cases}
\end{equation}
where $\mathbf{u}$ and $p$ are the velocity and pressure solutions, $\mathbf{f}$ is the external force, and $\mathbf{g}$ and $\mathbf{t}$, resp., are the Dirichlet and Neumann data, resp. (with $\int_{\partial\Omega} \mathbf{g}\cdot \mathbf{n} = 0$ if $\Gamma_D=\partial\Omega$). All physical quantities as dimensionless, with respect to a characteristic velocity $U$ and a characteristic length $L$.
The quantity $1/\nu>0$ is the Reynolds number $Re$. The pressure sign is chosen as in \cite{Navier-Stokes_1}. Note that $(\mathbf{\nabla}\mathbf{u})\mathbf{u}= (\mathbf{u}\cdot\mathbf{\nabla})\mathbf{u}$, where the tensor $\mathbf{\nabla}\mathbf{u}$ is the Jacobian matrix of $\mathbf{u}$.

We will develop our theoretical analysis under the simplifying assumption of homogeneous Dirichlet conditions on the whole of $\partial\Omega$ (whereas numerical results will be given in the general case). Therefore, we introduce the spaces
\begin{equation*}
\mathbf{V}:=[H^1_0(\Omega)]^2 \,, \qquad Q:=L^2_0(\Omega)=\big{\{}q\in L^2(\Omega) \quad\text{s.t.}\quad\int_{\Omega}q\, \text{d}\Omega=0\big{\}} \,,
\end{equation*}
with norms $ \|\mathbf{v}\|_\mathbf{V}= \Vert \mathbf{v} \Vert_{1,\Omega}$ and $\Vert q \Vert_Q = \Vert q \Vert_{0,\Omega}$, respectively, and we formulate the problem variationally as follows: {\em find $(\mathbf{u},p)\in \mathbf{V}\times Q$ such that}
\begin{equation}\label{NS_variational}
    \begin{cases}
    a(\mathbf{u},\mathbf{v})+c(\mathbf{u};\mathbf{u},\mathbf{v})+b(\mathbf{v},p)=(\mathbf{f},\mathbf{v}) \quad \forall \,\mathbf{v}\in \mathbf{V} \,, \\
    b(\mathbf{u},q)=0 \quad \forall\, q\in Q \,, 
    \end{cases}
\end{equation}
with 
\begin{equation}\label{eq:a_e_b}
\begin{split}
      &a(\mathbf{u},\mathbf{v}):=\int_{\Omega}\nu\,\mathbf{\nabla}\mathbf{u}:\mathbf{\nabla}\mathbf{v}\,\text{d}\Omega, \quad \forall\, \mathbf{u},\mathbf{v} \in \mathbf{V} \,, \\      
   &b(\mathbf{v},q):=\int_{\Omega} (\nabla \cdot \mathbf{v})\, q\, \text{d}\Omega \quad \forall \mathbf{v}\in\mathbf{V}, q\in Q \,, \\
   &c(\mathbf{w};\mathbf{u},\mathbf{v}):=\int_{\Omega}(\mathbf{\nabla}\mathbf{u})\mathbf{w}\,\cdot\,\mathbf{v}\,d\Omega \quad \forall \,\mathbf{w},\mathbf{u},\mathbf{v}\in\mathbf{V} \,, \\
   &  (\mathbf{f},\mathbf{v}):=\int_{\Omega}\mathbf{f}\cdot\mathbf{v}\,d\Omega \quad \forall \mathbf{v}\in\mathbf{V}, \text{ with } \mathbf{f}\in (L^2(\Omega))^2 \,.
\end{split}
\end{equation}
It is well known that the bilinear forms $a$ and $b$, as well as the trilinear form $c$, are continuous in their spaces of definition; furthermore, $a$ is coercive on $\mathbf{V}$ with coercivity constant $\nu$, whereas $b$ satisfies an inf-sup condition on $\mathbf{V}\times Q$ for some inf-sup constant $\beta>0$.
%The following properties hold:
%\begin{enumerate}
%    \item $a(\cdot,\cdot)$ is continuous and coercive:
%    \begin{equation}\label{eq:a_cont}
%        \lvert a(\mathbf{u},\mathbf{v})\rvert\leq C_a  \|\mathbf{u}\|_{\mathbf{V}}\|\mathbf{v}\|_{\mathbf{V}}\quad \forall\, \mathbf{u},\mathbf{v}\in \mathbf{V}
%    \end{equation}
%    \begin{equation}\label{eq:a_coer}
%        a(\mathbf{v},\mathbf{v})\geq \nu\|\mathbf{v}\|^2_{\mathbf{V}}\quad \forall q\in Q
%    \end{equation}
%    \item $b(\cdot,\cdot)$ is continuous and satisfies the inf-sup condition:
%    \begin{equation}\label{eq:b_cont}
%        \lvert b(\mathbf{v},q)\rvert\leq \|\mathbf{v}\|_{\mathbf{V}} \|q\|_{Q} \quad \forall \mathbf{v}\in\mathbf{V},\,q\in Q
%    \end{equation}
%    \begin{equation}\label{eq:b_inf_sup}
%       \exists \beta\quad\text{s.t.}\quad\sup_{\mathbf{v}\in\mathbf{V},\mathbf{v}\not=\mathbf{0}}\frac{b(\mathbf{v},q)}{\|\mathbf{v}\|_{\mathbf{V}}}\geq\beta \|q\|_Q\quad\forall q\in Q
%   \end{equation}
%\item The trilinear form $c(\cdot;\cdot,\cdot)$ is continuous i.e.:
%\begin{equation}\label{eq:C_cont}
%    c(\mathbf{w};\mathbf{u},\mathbf{v})\leq \hat{C} \|\mathbf{w}\|_{\mathbf{V}}\|\mathbf{u}\|_{\mathbf{V}}\|\mathbf{v}\|_{\mathbf{V}}
%\end{equation}
%\end{enumerate}

Throughout the paper, we will assume small data, precisely
\begin{equation}\label{gamma_small}
    \gamma:=\frac{\hat{C} \|\mathbf{f}\|_{-1,\Omega}}{\nu^2}<1 \,,
\end{equation}
where $\hat{C}$ is the continuity constant of the form $c$.
Under this assumption, the solution exists and is unique, moreover the velocity $\mathbf{u}$  satisfies
\begin{equation}\label{eq:impl_small_data}
    \|\mathbf{u}\|_\mathbf{V}\leq \frac{\|\mathbf{f}\|_{-1,\Omega}}{\nu} \,.
\end{equation}

\section{The discretization}\label{sec:discr}
We adopt the VEM discretization proposed in \cite{Navier-Stokes_1} (see also \cite{Stokes_1} and \cite{Darcy} for more details) with order of accuracy $k \geq 2$. We introduce a partition $\mathcal{T}_h$ made by polygonal elements $E$ of diameter $h_E$; we will set $h_\Omega= \max_{E \in \mathcal{T}_h} h_E$.   The elements $E$ satisfy two geometrical assumption which will be important in the proof related to the error estimator, even if we will restrict the numerical cases to a more specific class of elements.
The two hypotheses are:
\begin{itemize}\label{geom_ip}
    \item $E$ is star-shaped with respect to a ball $B_E$ of radius $\geq \eta h_E$\,;
    \item the distance between any two vertexes of $E$ is $\geq \sigma h_E$.
\end{itemize}
Each polygon $E\in\mathcal{T}_h$ has a certain number of edges $e\subset \partial E$, whose length will be denoted by $h_e$. $\mathcal{F}_h$ indicates the collection of all edges of the mesh, whereas $\mathcal{F}_h^{o}$ is the subset of the internal edges.

%Give a function $v$, square brackets $[v]$ represent the "jump" of the function $v$ on a edge i.e. the difference of the function values in the two elements sharing the edge considered. Given a vector-valued function $\mathbf{v}$, $[\mathbf{v}]\mathbf{n}_e=(\mathbf{v}_{E^+}-\mathbf{v}_{E^-})\mathbf{n}_e$ where $\mathbf{n}_e$ is the $e$-normal vector from $E^-$ to $E^+$.
%With $D_E$ we indicate the union of the element $E$ and the elements sharing at least one edge with $E$.

\subsection{VEM spaces}

At first, we define some important projections for the VEM formulation.

The \textit{Nabla-projection} $\mathbf{\Pi}^{\nabla,E}_k:\mathbf{V}\rightarrow[\mathbb{P}_k(E)]^2$ is defined by the conditions:
\begin{equation*}
    \begin{cases}
    \int_E \mathbf{\nabla}(\mathbf{v}-\mathbf{\Pi}^{\nabla,E}_k\mathbf{v}):\mathbf{\nabla}\mathbf{q} \, dE=0 \quad \forall \mathbf{v}\in\mathbf{V}\quad\text{and}\quad \forall \mathbf{q}\in[\mathbb{P}_k(E)]^2,\\
    \int_{\partial E}(\mathbf{v}-\mathbf{\Pi}_k^{\nabla,E}\mathbf{v}) =\mathbf{0} \,.
    \end{cases}
\end{equation*}

The \textit{$L^2$-projection} for scalar functions $\Pi_k^{0,E}:L^2(E)\rightarrow\mathbb{P}_k(E)$ ($k \geq 0$) is defined as:
\begin{equation*}
    \int_E (v-\Pi_k^{0,E}v)q \, dE=0\qquad  \forall v\in L^2(E)\quad\text{and}\quad\forall q\in\mathbb{P}_k(E) \,,
\end{equation*}
with natural extension to vector-valued case and tensor-valued case.

Following \cite{Stokes_1}, we recall the following definitions of spaces:
\begin{equation}\label{B_E}
    \mathbb{B}_k(E):=\{v\in C^0(\partial E) \quad\text{s.t.}\quad v_{\vert e}\in\mathbb{P}_k(e)\quad\forall e\subset \partial E\},
\end{equation}
\begin{equation*}
    \mathcal{G}_k(E):= \nabla (\mathbb{P}_{k+1}(E))\subseteq [\mathbb{P}_k(E)]^2 \,,
\end{equation*}
\begin{equation*}
    \mathcal{G}_k^{\oplus}(E):=\mathbf{x}^{\perp}[\mathbb{P}_{k-1}(E)]\subseteq [\mathbb{P}_k(E)]^2\quad\text{with}\quad\mathbf{x}^{\perp}:=(x_2,-x_1).
\end{equation*}

\smallskip
The local VEM space is defined starting from introducing
\begin{equation*}
\begin{split}
        \mathbf{U}_h^E=\{\mathbf{v}\in[H^1(E)]^2\quad\text{s.t.} \quad &\mathbf{v}_{\lvert\partial E}\in[\mathbb{B}_k(E)]^2\,, \quad 
        \nabla \cdot \mathbf{v} \in\mathbb{P}_{k-1}(E)\,, \\
        &-\mathbf{\Delta}\mathbf{v}-\nabla s \in \mathcal{G}^{\oplus}_k(E) \text{ for some } s\in L^2(E)  \ \} \,.
\end{split}        
\end{equation*}
Then, the VEM space for velocities associated with $E$ is
\begin{equation}\label{eq:VEM_space}
    \mathbf{V}^E_h=\{\mathbf{v}\in\mathbf{U}_h^E\quad\text{s.t.}\quad\big(\mathbf{v}-\mathbf{\Pi}_k^{\nabla,E}\mathbf{v},\mathbf{g}_k^{\perp}\big)_E=0 \quad \forall \mathbf{g}_k^{\perp}\in\mathcal{G}_k^{\oplus}(E)/\mathcal{G}_{k-2}^{\oplus}(E)     \} \,,
\end{equation}
where $\mathcal{G}_k^{\oplus}(E)/\mathcal{G}_{k-2}^{\oplus}(E) $ denotes the subspace of $\mathcal{G}_k^{\oplus}(E)$ which is $L^2(E)$-orthogonal to $\mathcal{G}_{k-2}^{\oplus}(E)$ (for $k=2$ it coincides with $\mathcal{G}_2^{\oplus}(E)$). Clearly $\mathbf{V}^E_h$ contains $[\mathbb{P}_2(E)]^2$ as a proper subspace.

According to \cite{Navier-Stokes_1}, one has
\begin{equation*}
    \text{dim}\mathbf{V}_h^E=2kn_E+k^2-k \,,
\end{equation*}
where $n_E$ denotes the number of vertices (or edges) of $E$.
It is possible to identify a function $\mathbf{v}\in\mathbf{V}^E_h$ by its degrees of freedom $(Dv)$, which can be divided into four groups (\cite{Stokes_1}), namely
\begin{enumerate}
    \item $Dv1$: the values of $\mathbf{v}$ at the vertices of $E$;
    \item $Dv2$: the values of $\mathbf{v}$ at $k-1$ interior points of each edge of $E$;
    \item $Dv3$: the following moments of $\mathbf{v}$:
    \begin{equation*}
        \int_E \mathbf{v}\, \mathbf{g}_{k-2}   \, dE \quad \forall  \mathbf{g}_{k-2} \in \mathcal{G}_{k-2}^{\oplus}(E) 
    \end{equation*}
    (this set of d.o.f. is empty when $k=2$);
    \item $Dv4$: the following moments of the divergence of $\mathbf{v}$ in $E$:
    \begin{equation*}
        \int_E (\nabla \cdot \mathbf{v})\, q_{k-1}\, dE \quad \forall q_{k-1} \in \mathbb{P}_{k-1}(E)/\mathbb{R} \,.
    \end{equation*}
\end{enumerate}
The space of pressures is chosen as
\begin{equation}\label{Stokes_Q_space}
    Q^E_h:=\mathbb{P}_{k-1}(E) \,,
\end{equation}
which has dimension $\tfrac12k(k+1)$; for any $q\in Q_h^E$ the degrees of freedom $(Dq)$ are the moments of  order up to $k-1$ of $q$ in $E$:
    \begin{equation*}
        \int_E q\, p_{k-1}\, dE\quad \forall p_{k-1}\in\mathbb{P}_{k-1}(E) \,.
    \end{equation*}

\subsection{Multilinear forms}
Multilinear forms are obtained by summing up elemental contributions defined on the discrete spaces introduced above, in such a way we can compute them. First of all we underline that the form $b$ is not approximated, since we can compute it element-wise using the local degrees of freedom.

The restiction $a^E $ of the form $a$ on the element $E\in\mathcal{T}_h$ is approximated  by $a_h^E: \mathbf{V}_h^E\times  \mathbf{V}_h^E\rightarrow \mathbb{R}$, defined as
\begin{equation*}
    a_h^E(\mathbf{u}_h,\mathbf{v}_h):=a^E(\mathbf{\Pi}_k^{\nabla,E}\mathbf{u}_h,\mathbf{\Pi}_k^{\nabla,E}\mathbf{v}_h)+{S}_h^E\big{(}(\mathbf{I}-\mathbf{\Pi}_k^{\nabla,E})\mathbf{u}_h,(\mathbf{I}-\mathbf{\Pi}_k^{\nabla,E})\mathbf{v}_h\big{)} \,,
\end{equation*}
where ${S}_h^E:\mathbf{V}_h^E\times\mathbf{V}_h^E\rightarrow\mathbb{R}$ is a stabilizing bilinear form satisfying
\begin{equation*}
    \gamma_* a^E(\mathbf{v}_h,\mathbf{v}_h)\leq {S}_h^E(\mathbf{v}_h,\mathbf{v}_h)\leq  \gamma^* a^E(\mathbf{v}_h,\mathbf{v}_h) \quad \forall\, \mathbf{v}_h\in\mathbf{V}_h\quad \text{s.t.}\quad\mathbf{\Pi}_k^{\nabla,E}\mathbf{v}_h=\mathbf{0} \,,
\end{equation*}
for some constants $0 < \gamma_* \leq \gamma^*$ independent of $h$ and $E$. Precisely, we choose as in \cite{Stokes_1}
\begin{equation*}
    S_h^E(\mathbf{u}_h,\mathbf{v}_h)=\alpha^E\bar{\mathbf{u}}_h\bar{\mathbf{v}}_h 
\end{equation*}
with $\alpha^E=\nu$, where $\bar{\mathbf{u}}_h, \bar{\mathbf{v}}_h$ are the vectors containing the values of the local degrees of freedom of $\mathbf{u}_h$ and $\mathbf{v}_h$. 

The form $a^E_h$ satisfies the {\em $k$-consistency} and {\em stability} properties:
\begin{equation*}     
\begin{split}
      &  a^E_h(\mathbf{v}_h,\mathbf{q}_k)=a^E(\mathbf{v}_h,\mathbf{q}_k) \qquad \forall \, \mathbf{v}_h\in\mathbf{V}_h^E\,, \quad \forall \mathbf{q}_k\in[\mathbb{P}_k(E)]^2  \,, \\
     &     \alpha_* a^E(\mathbf{v}_h,\mathbf{v}_h)\leq a_h^E(\mathbf{v}_h,\mathbf{v}_h)\leq  \alpha^* a^E(\mathbf{v}_h,\mathbf{v}_h) \quad \forall\, \mathbf{v}_h\in\mathbf{V}_h^E \,,
   \end{split}
    \end{equation*}
where $\alpha_*$ and $\alpha^*$ are positive constants independent of $h$ and $E$.

\medskip
The trilinear form $c$ is approximated on the element $E$ by the form 
$c_h^E: \mathbf{V}_h^E\times\mathbf{V}_h^E\times\mathbf{V}_h^E\rightarrow\mathbb{R}$ defined as in \cite{Navier-Stokes_1} as follows:  for any $\mathbf{w}_h,\mathbf{u}_h,\mathbf{v}_h\in\mathbf{V}_h^E$,
\begin{equation}\label{eq:c^E_def}
    c_h^E(\mathbf{w}_h;\mathbf{u}_h,\mathbf{v}_h):=\int_E\left(\left(\mathbf{\Pi}_{k-1}^{0,E}\mathbf{\nabla}\mathbf{u}_h\right)\left(\mathbf{\Pi}_k^{0,E}\mathbf{w}_h\right)\right)\cdot \mathbf{\Pi}_k^{0,E}\mathbf{v}_h\,dE \,.
    \end{equation}
Note that all projectors appearing in $c_h^E$ are computable, in particular it is possible to compute $\mathbf{\Pi}_k^{0,E}$ thanks to the choice of an enhanced space.

\medskip
Finally, the right-hand side $\mathbf{f}$ is locally approximated by $\mathbf{f}_h:=\mathbf{\Pi}^{0,E}_k \mathbf{f}$, which gives for any $\mathbf{v}_h\in\mathbf{V}_h$
\begin{equation}\label{rhs_NS}
    (\mathbf{f}_h,\mathbf{v}_h)=\sum_{E\in\mathcal{T}_h}\int_{E}\mathbf{f}_h\cdot\mathbf{v}_h\,dE=\sum_{E\in\mathcal{T}_h}\int_{E}\mathbf{\Pi}^{0,E}_k \mathbf{f}\cdot\mathbf{v}_h\,dE=\sum_{E\in\mathcal{T}_h}\int_{E}\mathbf{f}\cdot\mathbf{\Pi}^{0,E}_k\mathbf{v}_h\,dE \,.
\end{equation}
Since $\mathbf{\Pi}^{0,E}_k\mathbf{v}_h$ is computable, the integral can be computed by a quadrature formula.

\subsubsection{The discrete problem}
Now, it is possible to define the discrete problem: {\em find  $(\mathbf{u}_h,p_h)\in\mathbf{V}_h\times Q_h $ such that }
\begin{equation}\label{eq:NS_discr}
\begin{cases}
        a_h(\mathbf{u}_h,\mathbf{v}_h)+ {c}_h(\mathbf{u}_h;\mathbf{u}_h,\mathbf{v}_h)+b(\mathbf{v}_h,p_h)=(\mathbf{f}_h,\mathbf{v}_h)\quad\forall\,\mathbf{v}_h\in\mathbf{V}_h \,, \\
        b(\mathbf{u}_h,q_h)=0\quad\forall\,q_h\in Q_h \,.
    \end{cases}
\end{equation}
It is worth noticing that the discrete velocity $\mathbf{u}_h$ is exactly {\em divergence-free} in $\Omega$, as a consequence of the identity $\text{div}\mathbf{V}_h=Q_h$ (see \cite{Navier-Stokes_1}).

\begin{remark}
The mathematical analysis of the discrete problem is facilitated if the convective term is approximated by the {\em skew-symmetric} form $\tilde{c}_h(\mathbf{w}_h;\mathbf{u}_h,\mathbf{v}_h):=\frac12 {c}_h(\mathbf{w}_h;\mathbf{u}_h,\mathbf{v}_h) - \frac12 {c}_h(\mathbf{w}_h;\mathbf{v}_h,\mathbf{u}_h)$. In this case, existence and uniqueness follows from the assumption of data smallness, as indicated by the following result. 
\begin{theorem}[\cite{Navier-Stokes_1},{Theorem 3.5}]
If 
\begin{equation}\label{eq:as_gamma_h_NS}
    \gamma_h:=\frac{\hat{C}_h\|\mathbf{f}_h\|_{-1,\Omega}}{\alpha_*^2\nu^2}\leq r <1\,,
\end{equation}
the VEM discretization \eqref{eq:NS_discr} with $c_h$ replaced by $\tilde{c}_h$ has a unique solution $(\mathbf{u}_h,p_h)\in\mathbf{V}_h\times Q_h$, such that 
\begin{equation}\label{eq:u_err}
    \|\mathbf{u}_h\|_{\mathbf{V}}\leq\frac{\|\mathbf{f}_h\|_{-1,\Omega}}{\alpha_*\nu} \,.
\end{equation}
\end{theorem}

Furthermore, if $\mathbf{u},\mathbf{f}\in[H^{s+1}(\Omega)]^2$ and $p \in H^s(\Omega)$ for some $0\leq s\leq k$, then the following a priori error bounds hold:
\begin{equation}\label{eq:apriori-bound}
\begin{split}
&   \|\mathbf{u}-\mathbf{u}_h\|_\mathbf{V}\leq h^s \mathcal{N}(\mathbf{u};\nu,\gamma,r)+ h^{s+2}\mathcal{S}(\mathbf{f};\nu,r) \,, \\
&  \|p-p_h\|_Q\leq C h^s \lvert p\rvert_{s,\Omega}+C h^{s+2} \lvert\mathbf{f}\rvert_{s+1,\Omega}+h^s\mathcal{H}(\mathbf{u};\nu,\gamma,r) \,,
\end{split}
\end{equation}
for suitable constants $\mathcal{N}$, $\mathcal{S}$, $\mathcal{H}$ independent of $h$.

\medskip
Results of this type can be proven also for the discretization \eqref{eq:NS_discr} which does not use the skew-symmetric form of the convective term (see \cite[Remark 3.5]{Navier-Stokes_1}), at the expense of a more involved analysis. 
\end{remark}

\section{A posteriori error analysis}
\label{sec:main}
Let us set
\begin{equation}
    \mathbf{e}_{\mathbf{u}}=\mathbf{u}-\mathbf{u}_h\,,  \qquad e_p=p-p_h\,,
\end{equation}
where $(\mathbf{u},p)$ is the solution of \eqref{NS_variational}
and $(\mathbf{u}_h,p_h)$ is the corresponding solution of \eqref{eq:NS_discr}. In order to build a proper a posteriori error estimator, we partially follow the approach of \cite{Ains}. 
\subsection{Residual analysis}
For any $\mathbf{v}\in\mathbf{V}$ and $q\in Q$, let us define
\begin{equation}\label{eq:R_NS}
    \mathcal{R}^{NS}[(\mathbf{e}_{\mathbf{u}},e_p),(\mathbf{v},q)]= a(\mathbf{e}_{\mathbf{u}},\mathbf{v})+b(\mathbf{v},e_p)+b(\mathbf{e}_{\mathbf{u}},q)+\delta(\mathbf{u},\mathbf{u}_h,\mathbf{v}) \,,
\end{equation}
where
\begin{equation}
    \delta(\mathbf{u},\mathbf{u}_h,\mathbf{v})=c(\mathbf{u},\mathbf{u},\mathbf{v})-c(\mathbf{u}_h,\mathbf{u}_h,\mathbf{v}) \,.
\end{equation}
It is useful to set
\begin{equation}
    \mathcal{L}^{NS}[(\mathbf{u},p),(\mathbf{v},q)]=a({\mathbf{u}},\mathbf{v})+b(\mathbf{v},p)+b(\mathbf{u},q)+c(\mathbf{u},\mathbf{u},\mathbf{v}) \quad \forall (\mathbf{v},q)\in \mathbf{V}\times Q
\end{equation}
and
\begin{equation}
\begin{split}
    \mathcal{L}_h^{NS}[(\mathbf{u}_h,p_h),(\mathbf{v}_h,q_h)]=a_h({\mathbf{u}}_h,\mathbf{v}_h)&+b(\mathbf{v}_h,p_h) + b(\mathbf{u}_h,q_h) +\\
    &+ c_h(\mathbf{u}_h,\mathbf{u}_h,\mathbf{v}_h) \quad \forall (\mathbf{v}_h,q_h)\in \mathbf{V}_h\times Q_h \,.
\end{split}
\end{equation}
Following \cite{Ains}, let us introduce the pair $(\mathbf{\Phi},\psi)\in \mathbf{V}\times Q$ defined by
\begin{equation}\label{eq:a_g_R}
    a(\mathbf{\Phi},\mathbf{v})+g(\psi,q)=\mathcal{R}^{NS}[(\mathbf{e}_{\mathbf{u}},e_p),(\mathbf{v},q)] \quad \forall (\mathbf{v},q)\in \mathbf{V}\times Q \,,
\end{equation}
where $g: Q\times Q \rightarrow \mathbb{R}$, $g(p,q)=\int_{\Omega}p\,q\,d\Omega$.
The existence and uniqueness of $(\mathbf{\Phi},\psi)$ 
is guaranteed since $\mathcal{R}^{NS}[(\mathbf{e}_{\mathbf{u}},e_p), \ \cdot \ ]$ defines a linear continuous functional on $\mathbf{V}\times Q$.
Now it is convenient to define the norm
\begin{equation}
    \vert \! \vert \! \vert(\mathbf{\Phi},\psi) \vert \! \vert \! \vert =(a(\mathbf{\Phi},\mathbf{\Phi})+g(\psi,\psi))^{1/2} =(\nu\vert\mathbf{\Phi}\vert^2_{1,\Omega}+\|\psi\|_{0,\Omega}^2)^{1/2} \,.
\end{equation}
It is shown in \cite{Ains} that this norm and the norm $(\nu\vert\mathbf{e}_{\mathbf{u}}\vert_{1,\Omega}^2+\|e_p\|_{0,\Omega}^2)^{1/2}$ are equivalent. 
\begin{theorem}[Norm equivalence]\label{thm:norm_eq}
  There exist two positive constants $C_1$ and $C_2$ such that
  \begin{equation}\label{eq:norm_eq}
      C_1 \vert \! \vert \! \vert(\mathbf{\Phi},\psi) \vert \! \vert \! \vert^2\leq \nu\vert\mathbf{e}_\mathbf{u}\vert_{1,\Omega}^2+\|e_p\|_{0,\Omega}^2\leq C_2 \vert \! \vert \! \vert(\mathbf{\Phi},\psi) \vert \! \vert \! \vert^2 \,.
  \end{equation}
 The constant $C_2$ depends on $\hat{C}$, $\beta$, $\gamma$, whereas $C_1$ depends on $\hat{C}$ and $\nu$.
\end{theorem}
%\vskip -.5cm
\noindent It is important to notice that the proof uses the hypothesis \eqref{gamma_small} of having small data, and a convergence result of the type \eqref{eq:apriori-bound}.

In our situation, we actually have $\psi=0$, as a consequence of the property that both $\mathbf{u}$ and $\mathbf{u}_h$ are divergence-free. Indeed, this implies $b(\mathbf{e_u},q)=0\quad\forall q\in Q$, and taking $(\mathbf{0},q)$ as test function in \eqref{eq:a_g_R} yields $g(\psi,q)=0 \quad \forall q \in Q$, i.e., $\psi=0$. Combining this with \eqref{eq:norm_eq} and the identity
\begin{equation}\label{eq:phi-identity}
\vert \! \vert \! \vert(\mathbf{\Phi},0) \vert \! \vert \! \vert^2 = a(\mathbf{\Phi},\mathbf{\Phi})= \mathcal{R}^{NS}[(\mathbf{e}_{\mathbf{u}},e_p),(\mathbf{\Phi},0)]
\end{equation}
which follows from \eqref{eq:a_g_R} with $(\mathbf{v},q)=(\mathbf{\Phi},0)$, we arrive at the following error bound:
\begin{equation}\label{eq:errbnd-1}
\nu\vert\mathbf{e}_\mathbf{u}\vert_{1,\Omega}^2+\|e_p\|_{0,\Omega}^2\leq C_2 \mathcal{R}^{NS}[(\mathbf{e}_{\mathbf{u}},e_p),(\mathbf{\Phi},0)] \,.
\end{equation}
Thus, we are left with the task of bounding the right-hand side.

%\medskip
To proceed, we take inspiration from the analogous Stokes case treated in \cite{a_post_Stokes}, namely we rewrite the right-hand side as follows:
\begin{equation*}
\begin{split}
     \mathcal{R}^{NS}[(\mathbf{e}_{\mathbf{u}},e_p),(\mathbf{\Phi},0)] &=(\mathbf{f},\mathbf{\Phi})-\mathcal{L}^{NS}[(\mathbf{u}_h,p_h),(\mathbf{\Phi},0)]
     \\
     &=(\mathbf{f},\mathbf{\Phi})-(\mathbf{f}_h,\boldsymbol{\xi}_h) +\mathcal{L}_h^{NS}[(\mathbf{u}_h, p_h), (\boldsymbol{\xi}_h,0)] \\
     & \qquad  -\mathcal{L}^{NS}[(\mathbf{u}_h, p_h), (\boldsymbol{\xi}_h,0)] -\mathcal{L}^{NS}[(\mathbf{u}_h, p_h), (\boldsymbol{\Phi}-\boldsymbol{\xi}_h,0)]\\
     &=(\mathbf{f}-\mathbf{f}_h,\mathbf{\Phi})+(\mathbf{f}_h,\mathbf{\Phi}-\boldsymbol{\xi}_h) \\
     & \qquad +\mathcal{L}_h^{NS}[(\mathbf{u}_h, p_h), (\boldsymbol{\xi}_h,0)] -\mathcal{L}^{NS}[(\mathbf{u}_h, p_h), (\boldsymbol{\xi}_h,0)] \\
     &\qquad  -\mathcal{L}^{NS}[(\mathbf{u}_h, p_h), (\boldsymbol{\Phi}-\boldsymbol{\xi}_h,0)]  \,,
\end{split}
\end{equation*}
where $\boldsymbol{\xi}_h$ is a generic function belonging to $\mathbf{V}_h$.
Focusing on the last terms, we notice that
\begin{equation}\label{eq:L-L_h}
\begin{split}
   \mathcal{L}_h^{NS}[(\mathbf{u}_h, p_h), (\boldsymbol{\xi}_h,0)]-\mathcal{L}^{NS}[(\mathbf{u}_h, p_h), (\boldsymbol{\xi}_h,0)]=a_h(\mathbf{u}_h,\boldsymbol{\xi}_h)-a(\mathbf{u}_h,\boldsymbol{\xi}_h)\\
+ c_h(\mathbf{u}_h;\mathbf{u}_h,\boldsymbol{\xi}_h)-c(\mathbf{u}_h;\mathbf{u}_h,\boldsymbol{\xi}_h) \,;
\end{split}
\end{equation} 
introducing, for any edge $e\in\mathcal{F}_h^o$ which is shared by two elements $E_1$ and $E_2$ with outward unit normal vectors $\mathbf{n}_e^{E_1}$ and $\mathbf{n}_e^{E_2}$, the normal jump 
$[\![ \mathbf{v}]\!]_e = \mathbf{v}_{\vert E_1} \mathbf{n}_e^{E_1} + \mathbf{v}_{\vert E_2} \mathbf{n}_e^{E_2}$, and
setting for convenience $\mathbf{w}=\boldsymbol{\Phi}-\boldsymbol{\xi}_h$, we have

\begin{equation}\label{eq_L_NS}
    \begin{split}
        &\mathcal{L}^{NS}[(\mathbf{u}_h, p_h), (\mathbf{w},0)]=a(\mathbf{u}_h,\mathbf{w})+b(\mathbf{w},p_h)+c(\mathbf{u}_h;\mathbf{u}_h,\mathbf{w})\\
        & \quad =\underbrace{\sum_{E\in\mathcal{T}_h}(-\nu\Delta\mathbf{\Pi}_k^{\nabla,E}\mathbf{u}_h-\nabla p_h,\mathbf{w})_E}_{(I)}
          +\underbrace{\sum_{e\in\mathcal{F}_h^o}([\![\nu\nabla \mathbf{\Pi}_k^{\nabla }\mathbf{u}_h+p_h\mathbf{I}]\!]_e,\mathbf{w})_e}_{(II)}\\
        &\qquad+\underbrace{\sum_{E\in\mathcal{T}_h}\nu(\nabla(\mathbf{u}_h-\mathbf{\Pi}_k^{\nabla,E}\mathbf{u}_h),\nabla\mathbf{w})_E}_{(III)}+ c(\mathbf{u}_h;\mathbf{u}_h,\mathbf{w})\\
        &\quad =(I)+(II)+(III)+\underbrace{c(\mathbf{u}_h;\mathbf{u}_h,\mathbf{w})-c_h(\mathbf{u}_h;\mathbf{u}_h,\mathbf{w})}_{(IV)}
           +c_h(\mathbf{u}_h;\mathbf{u}_h,\mathbf{w})\\
        &\quad =(I)+(II)+(III)+(IV)
        +\underbrace{\sum_{E\in\mathcal{T}_h} \big( (\mathbf{\Pi}^{0,E}_{k-1}\nabla\mathbf{u}_h ) \mathbf{\Pi}^{0,E}_k\mathbf{u}_h,\mathbf{\Pi}^{0,E}_k\mathbf{w}-\mathbf{w} \big)_E}_{(V)}\\
        &\qquad+\sum_{E\in\mathcal{T}_h} \big( (\mathbf{\Pi}^{0,E}_{k-1}\nabla\mathbf{u}_h ) \mathbf{\Pi}^{0,E}_k\mathbf{u}_h,\mathbf{w}\big)_E
        \\
        &\quad =\sum_{E\in\mathcal{T}_h} \big(-\nu\Delta\mathbf{\Pi}_k^{\nabla,E}\mathbf{u}_h-\nabla p_h+ (\mathbf{\Pi}^{0,E}_{k-1}\nabla\mathbf{u}_h ) \mathbf{\Pi}_k^{0,E}\mathbf{u}_h,\mathbf{w}\big)_E\\
        & \qquad +\sum_{e\in\mathcal{F}_h^o}([\![ \nu \nabla \mathbf{\Pi}_k^{\nabla}\mathbf{u}_h+p_h\mathbf{I}]\!]_e,\mathbf{w})_e
        +(III)+(IV)+(V).
    \end{split}
\end{equation}
Combining the term $(IV)$ in \eqref{eq_L_NS} with \eqref{eq:L-L_h}, we arrive at the following result.
\begin{lemma}  The right-hand side of \eqref{eq:errbnd-1} can be split as 
\begin{equation}\label{eq:RJi}
\begin{split}
     \mathcal{R}^{NS}[(\mathbf{e}_{\mathbf{u}},e_p),(\mathbf{\Phi},0)]=\sum_{i=1}^7 J_i \,,
\end{split}
\end{equation}
where
\begin{equation*}
    J_1=(\mathbf{f}-\mathbf{f}_h,\mathbf{\Phi}) \,,
\end{equation*}
\begin{equation*}
    J_2=a_h(\mathbf{u}_h,\boldsymbol{\xi}_h)-a(\mathbf{u}_h,\boldsymbol{\xi}_h) \,,
\end{equation*}
\begin{equation*}
    J_3=\sum_{E\in\mathcal{T}_h}\big(\mathbf{f}_h+\nu\Delta\mathbf{\Pi}^{\nabla,E}_k\mathbf{u}_h+\nabla p_h-(\mathbf{\Pi}^{0,E}_{k-1}\nabla\mathbf{u}_h) \mathbf{\Pi}^{0,E}_k\mathbf{u}_h,\mathbf{\Phi}-\boldsymbol{\xi}_h \big)_E \,,
\end{equation*}
\begin{equation*}
    J_4=-\sum_{e\in\mathcal{F}_h^o}([\![\nu\nabla \mathbf{\Pi}^{\nabla}_k\mathbf{u}_h+p_h\mathbf{I}]\!]_e,\mathbf{\Phi}-\boldsymbol{\xi}_h)_e \,,
\end{equation*}
\begin{equation*}
    J_5=-\sum_{E\in\mathcal{T}_h}\nu \big(\nabla(\mathbf{u}_h-\mathbf{\Pi}^{\nabla,E}_k\mathbf{u}_h),\nabla(\mathbf{\Phi}-\boldsymbol{\xi}_h) \big)_E \,,
\end{equation*}
\begin{equation*}
    J_6=c_h(\mathbf{u}_h;\mathbf{u}_h,\mathbf{\Phi})-c(\mathbf{u}_h;\mathbf{u}_h,\mathbf{\Phi}) \,,
\end{equation*}
\begin{equation*}
    J_7=\sum_{E\in\mathcal{T}_h}\big( (\mathbf{\Pi}^{0,E}_{k-1}\nabla\mathbf{u}_h) \mathbf{\Pi}^{0,E}_k\mathbf{u}_h,(\mathbf{\Phi}-\boldsymbol{\xi}_h)-\mathbf{\Pi}^{0,E}_k(\mathbf{\Phi}-\boldsymbol{\xi}_h)\big)_E \,.
\end{equation*}
\end{lemma}
Comparing the analogous decomposition done for the Stokes case (see \cite{a_post_Stokes}), we observe that the term $J_3$ changes showing the additional convective term, whereas  two new terms, $J_6$ and $J_7$, appear.

\subsection{The a posteriori estimator}\label{sec:estimator}

We are ready to introduce our element-wise estimator, which is defined additively as follows:
\begin{equation}\label{eq:def-components}
\eta_E^2=\eta_{f,E}^2+\eta_{B,E}^2+\eta_{e,E}^2+\eta_{S,E}^2+\eta_{c1,E}^2+\eta_{c2,E}^2+\eta_{c3,E}^2\,,
\end{equation}
where
\begin{equation*}
\eta_{f,E}^2=h_E^2 \|\mathbf{f}_h-\mathbf{f}\|_{0,E}^2 \,; 
\end{equation*}
\begin{equation*}
\eta_{B,E}^2=h_E^2 \Vert\mathbf{f}_h+\nu\Delta\mathbf{\Pi}^{\nabla,E}_k\mathbf{u}_h+\nabla p_h- (\mathbf{\Pi}^{0,E}_{k-1}\nabla\mathbf{u}_h) \mathbf{\Pi}^{0,E}_k\mathbf{u}_h \Vert_{0,E}^2 \,;
\end{equation*}
\begin{equation*}
\eta_{e,E}^2=\sum_{e\in E}h_e\big\Vert[\![\nu\nabla \mathbf{\Pi}^{\nabla}_k\mathbf{u}_h+p_h\mathbf{I}]\!]_e\big\Vert^2_{0,e} \,;
\end{equation*}
\begin{equation*}
\eta_{S,E}^2=(\nu \sigma_E)^2
\qquad \text{with \ } 
\sigma_E=\big(\nu^{-1}S^E((\mathbf{I}-\mathbf{\Pi}^{\nabla,E}_k)\mathbf{u}_h,(\mathbf{I}-\mathbf{\Pi}^{\nabla,E}_k)\mathbf{u}_h)\big)^{1/2} \,;
\end{equation*}
\smallskip
\begin{equation*}
\eta_{c1,E}^2=h_E^2 \Vert \, \mathbf{\Pi}_k^{0,E}\big( (\mathbf{\Pi}_{k-1}^{0,E}\nabla\mathbf{u}_h) \mathbf{\Pi}_k^{0,E}\mathbf{u}_h\big)-(\mathbf{\Pi}_{k-1}^{0,E}\nabla\mathbf{u}_h) \mathbf{\Pi}_k^{0,E}\mathbf{u}_h  \Vert^2_{0,E}\,;
\end{equation*}
\begin{equation*}
\eta_{c2,E}^2=\sigma_E^2 \, 
    \Vert\mathbf{\Pi}_k^{0,E}\mathbf{u}_h \Vert^2_{L^\infty(E)}\,;
\end{equation*}
%where, for $\mathbf{w}\in[L^\infty(E)]^2$, we set $ \Vert \mathbf{w}\Vert_{\infty,E}=\max_{j\in\{1,2\}}\Vert {w}_j \Vert_{\infty,E} $;
\begin{equation*}
\begin{split}
\eta_{c3,E}^2=
     \Big(\sigma_E+\Vert\nabla(\mathbf{\Pi}_k^{0,E}\mathbf{u}_h-\mathbf{\Pi}_k^{\nabla,E}\mathbf{u}_h)\Vert_{0,E}\Big)^2
\left(\sigma_E+\Vert\nabla\mathbf{\Pi}^{\nabla,E}_k\mathbf{u}_h \Vert_{0,E} \right)^2\,.
\end{split}
\end{equation*}

\subsection{A posteriori upper bound of the error}\label{sec:estimate}

The following theorem guarantees the reliability of our error estimator \eqref{eq:def-components} in the energy norm.

\begin{theorem}[A posteriori error estimate]\label{theo:reliability}
Let $(\mathbf{u},p)\in\mathbf{V}\times Q$ be the solution of problem \eqref{NS_variational} and $(\mathbf{u}_h,p_h)\in\mathbf{V}_h\times Q_h$ be the solution of problem \eqref{eq:NS_discr}. The following a posteriori error estimate from above holds true:
\begin{equation}\label{eq:est_post}
    \nu\vert\mathbf{u}-\mathbf{u}_h\vert^2_{1,\Omega}+\|p-p_h\|_{0,\Omega}^2 \lesssim \frac{1}{\nu}\eta^2 \,,
\end{equation}
where
\begin{equation}
    \eta^2=\sum_{E\in\mathcal{T}_h}\eta^2_E \,.
\end{equation}
\end{theorem}
\begin{proof}
Recalling the bound \eqref{eq:errbnd-1} and the expression \eqref{eq:RJi}, we have to estimate each term $J_i$, $i=1, \dots, 7$.

Similarly to what was done for the Stokes case in \cite{a_post_Stokes},  we choose $\boldsymbol{\xi}_h=\mathbf{\Phi}_I$, the Scott-Zhang quasi-interpolant of $\mathbf{\Phi}$, which satisfies $\|\mathbf{\Phi}-\mathbf{\Phi}_I\|_{0,E}+h_E\vert\mathbf{\Phi}-\mathbf{\Phi}_I\vert_{1,E}\leq h_E\vert\mathbf{\Phi}\vert_{1,D_E}$ for all $E \in \mathcal{T}_h$, where $D_E$ denotes the union of $E$ and the elements that share at least one edge with $E$.

For the terms $J_2,J_4$ and $J_5$ we can use the results of the corresponding proof given in \cite{a_post_Stokes} for the analogous Stokes case. This can be done for $J_1$ as well, even if we are using a different expression of $\mathbf{f}_h$.
Repeating some calculations done in \cite{a_post_Stokes}, we obtain
\begin{equation}\label{eq:BD-J1}
\begin{split}
& J_1\lesssim\sum_{E\in\mathcal{T}_h}h_E\|\mathbf{f}-\mathbf{f}_h\|_E\vert\mathbf{\Phi}\vert_{1,E}\lesssim \big(\sum_{E\in\mathcal{T}_h}\eta^2_{f,E}\big)^{1/2}\vert\mathbf{\Phi}\vert_{1,\Omega} \,; \\
& J_2\lesssim\sum_{E\in\mathcal{T}_h}\nu \sigma_E\vert\mathbf{\Phi}\vert_{1,D_E}\lesssim \big(\sum_{E\in\mathcal{T}_h}\eta^2_{S,E}\big)^{1/2}\vert\mathbf{\Phi}\vert_{1,\Omega} \,; \\
 &   J_4 \lesssim \big(\sum_{e\in\mathcal{F}_h^o} h_e \|[\![\nu\nabla\mathbf{\Pi}^{\nabla}_k\mathbf{u}_h-p_h\mathbf{I}]\!]_e\|^2_{0,e}\big)^{\frac{1}{2}}\big(\sum_{E\in\mathcal{T}_h}\vert\mathbf{\Phi}\vert^2_{1,D_E}\big)^{\frac{1}{2}}\\
 & \qquad \lesssim \big(\sum_{E\in\mathcal{T}_h}\eta^2_{e,E}\big)^{1/2}\vert\mathbf{\Phi}\vert_{1,\Omega} \,; \\
& J_5\lesssim \sum_{E\in\mathcal{T}_h}\nu \sigma_E\vert\mathbf{\Phi}\vert_{1,D_E}\lesssim \big(\sum_{E\in\mathcal{T}_h}\eta^2_{S,E}\big)^{1/2}\vert\mathbf{\Phi}\vert_{1,\Omega} \,.
\end{split}
\end{equation}

% In fact using the definition of $\mathbf{f}_h$:
% \begin{equation}
% J_1=\sum_{E\in\mathcal{T}_h}(\mathbf{f}-\mathbf{f}_h,\mathbf{\Phi}-\mathbf{\Pi}_0^{0,E}\mathbf{\Phi})_E\lesssim\sum_{E\in\mathcal{T}_h}h_E\|\mathbf{f}-\mathbf{f}_h\|_E|\mathbf{\Phi}|_{1,E}
% \end{equation}
The term $J_3$ is upper bounded coherently with its new definition which includes the convection part:
\begin{equation}\label{eq:BD-J3}
\begin{split}
J_3&\lesssim\sum_{E\in\mathcal{T}_h} h_E \Vert\mathbf{f}_h+\nu\Delta\mathbf{\Pi}^{\nabla,E}_k\mathbf{u}_h+\nabla p_h-(\mathbf{\Pi}^{0,E}_{k-1}\nabla\mathbf{u}_h) \mathbf{\Pi}^{0,E}_k\mathbf{u}_h \Vert_{0,E}\vert\mathbf{\Phi}\vert_{1,D_E}\\
&\lesssim \big(\sum_{E\in\mathcal{T}_h}\eta^2_{B,E}\big)^{1/2}\vert\mathbf{\Phi}\vert_{1,\Omega} \,.
\end{split}
\end{equation}

Let us focus now on the terms $J_6$ and $J_7$. Let the trilinear form $c^*:[L^{4}(\Omega)]^2\times[L^{2}(\Omega)]^{2\times 2}\times [L^{4}(\Omega)]^2\rightarrow \mathbb{R}$ be defined as
\begin{equation}
    c^*(\mathbf{a};\mathbf{B},\mathbf{c})=\int_{\Omega}(\mathbf{B}\mathbf{a})\cdot \mathbf{c} \,d\Omega \,,
\end{equation}
with the obvious restriction $c^{*,E}$ on an element $E$. Moreover let us set
\begin{equation}
\mathbf{t}:=\left(\nabla\mathbf{u}_h\right)\mathbf{u}_h\,, \qquad     \mathbf{t}^E_h:=\left(\mathbf{\Pi}_{k-1}^{0,E}\nabla\mathbf{u}_h\right)\mathbf{\Pi}_k^{0,E}\mathbf{u}_h \,.
\end{equation}
Using the definition of $L^2$-projection $\mathbf{\Pi}_k^{0,E}$ upon $[\mathbb{P}_2(E)]^2$, we write $J_6$ as follows:
\begin{equation}\label{eq:split_j7}
\begin{split}
J_6 &= c_h(\mathbf{u}_h;\mathbf{u}_h,\mathbf{\Phi})-c(\mathbf{u}_h;\mathbf{u}_h,\mathbf{\Phi})=\sum_{E\in\mathcal{T}_h}
\left( \int_E \mathbf{t}^E_h\cdot\mathbf{\Pi}_k^{0,E}\mathbf{\Phi}\, d\,E-\int_E \mathbf{t}\cdot\mathbf{\Phi} \, d\,E \right)\\
&=\sum_{E\in\mathcal{T}_h} \left( \int_E \mathbf{\Pi}_k^{0,E}\mathbf{t}_h^E \cdot\mathbf{\Phi}\, d\,E-\int_E \mathbf{t}\cdot\mathbf{\Phi}\, dE \right)
=\sum_{E\in\mathcal{T}_h} (\mathbf{\Pi}_k^{0,E}\mathbf{t}_h^E-\mathbf{t},\mathbf{\Phi})_E \,.
\end{split}
\end{equation}
So,
\begin{equation}\label{eq:bound-J6}
\begin{split}
J_6&=\sum_{E\in\mathcal{T}_h} (\mathbf{\Pi}_k^{0,E}\mathbf{t}^E_h-\mathbf{t}^E_h+\mathbf{t}^E_h-\mathbf{t},\mathbf{\Phi})_E\\
&=\underbrace{\sum_{E\in\mathcal{T}_h} (\mathbf{\Pi}_k^{0,E}\mathbf{t}^E_h-\mathbf{t}^E_h,\mathbf{\Phi}-\mathbf{\Pi}_0^{0,E}\mathbf{\Phi})_E}_{(T_1)}+\sum_{E\in\mathcal{T}_h}(\mathbf{t}^E_h-\mathbf{t},\mathbf{\Phi})_E\\
  & =(T_1)+c^*(\mathbf{\Pi}^{0,E}_k\mathbf{u}_h;\mathbf{\Pi}^{0,E}_{k-1}\nabla \mathbf{u}_h,\mathbf{\Phi})-c(\mathbf{u}_h;\mathbf{u}_h,\mathbf{\Phi})\\  
& =(T_1)+ \sum_{E\in\mathcal{T}_h}c^{*,E}(\mathbf{\Pi}^{0,E}_k\mathbf{u}_h;\mathbf{\Pi}^{0,E}_{k-1}\nabla \mathbf{u}_h,\mathbf{\Phi})-\sum_{E\in\mathcal{T}_h}c^{*,E}(\mathbf{\Pi}^{0,E}_k\mathbf{u}_h;\nabla\mathbf{u}_h,\mathbf{\Phi}) +\\
&\quad+ \sum_{E\in\mathcal{T}_h}c^{*,E}(\mathbf{\Pi}^{0,E}_k\mathbf{u}_h;\nabla\mathbf{u}_h,\mathbf{\Phi}) - c(\mathbf{u}_h;\mathbf{u}_h,\mathbf{\Phi}) \\
&= (T_1)+\underbrace{\sum_{E\in\mathcal{T}_h}c^{*,E}(\mathbf{\Pi}^{0,E}_k\mathbf{u}_h;\mathbf{\Pi}_{k-1}^{0,E}\nabla\mathbf{u}_h-\nabla\mathbf{u}_h,\mathbf{\Phi})}_{(T_2)}\\
&\quad\quad\quad+\underbrace{\sum_{E\in\mathcal{T}_h}c^{*,E}(\mathbf{\Pi}^{0,E}_k\mathbf{u}_h-\mathbf{u}_h;\nabla\mathbf{u}_h,\mathbf{\Phi}}_{(T_3)}) \,.
\end{split}
\end{equation}
The term $(T_1)$ is easily bounded as follows:
\begin{equation}\label{eq_box1}
\begin{split}
    (T_1)&\leq \sum_{E\in\mathcal{T}_h} \Vert\mathbf{\Pi}_k^{0,E}\mathbf{t}^E_h-\mathbf{t}^E_h \Vert_{0,E} \Vert\mathbf{\Phi}-\mathbf{\Pi}_0^{0,E}\mathbf{\Phi}\Vert_{0,E}\\
    &\lesssim \sum_{E\in\mathcal{T}^E_h}\Vert\mathbf{\Pi}_k^{0,E}\mathbf{t}^E_h-\mathbf{t}^E_h \Vert_{0,E}h_E\vert\mathbf{\Phi}\vert_{1,E} 
        \lesssim \Big(\sum_{E\in\mathcal{T}_h}\eta^2_{c1,E}\Big)^{1/2}\vert\mathbf{\Phi}\vert_{1,\Omega} \,.
\end{split}
\end{equation}
Concerning  the term $(T_2)$, we start by using H\"old{e}r's inequality with exponents $\left(\frac{1}{2},\frac{1}{4},\frac{1}{4}\right)$ to get
\begin{equation*}
 (T_2) \lesssim  \sum_{i,j=1}^2\sum_{E\in\mathcal{T}_h}\Big\Vert\Pi_{1}^{0,E}\frac{\partial (\mathbf{u}_h)_i}{\partial x_j}-\frac{\partial (\mathbf{u}_h)_i}{\partial x_j}\Big\Vert_{0,E} \Vert\Pi_2^{0,E}(\mathbf{u}_h)_j\Vert_{L^4(E)} \Vert\mathbf{\Phi}_i \Vert_{L^4(E)} \,.
\end{equation*}
The Sobolev embedding $H^1(E)\subset L^4(E)$ with scaled continuous inclusion
$\|\cdot\|_{L^4(E)}\lesssim h_E^{-1/2}\|\cdot\|_{H^1(E)}$, together with the assumptions on the partition ${\cal T}_h$, yield
\begin{equation*}
\begin{split}
     (T_2) 
    &\lesssim \sum_{i,j=1}^2\Big(\sum_{E\in\mathcal{T}_h}\Big(\Big\Vert(\Pi_{1}^{0,E}\frac{\partial (\mathbf{u}_h)_i}{\partial x_j}-\frac{\partial (\mathbf{u}_h)_i}{\partial x_j}\Big\Vert_{0,E}\|\Pi_2^{0,E}(\mathbf{u}_h)_j\|_{L^4(E)}h_E^{-\frac{1}{2}}\Big)^2\Big)^{\frac{1}{2}}\times\\
    & \quad\times\Big(\sum_{E\in\mathcal{T}_h}\|\mathbf{\Phi}_i\|^2_{H^1(E)}\Big)^{\frac{1}{2}} \,.
    \end{split}
    \end{equation*}
On the other hand, the Lebesgue embedding $L^{\infty}(E)\subset L^4(E)$ with scaled continuous inclusion
$\|\cdot\|_{L^4(E)}\leq C h_E^{1/2} \|\cdot\|_{L^{\infty}(E)}$, together with the Poincar\'e inequality in $\Omega$, yield
\begin{equation*}
\begin{split}
     (T_2) & \lesssim \sum_{i,j=1}^2\Big(\sum_{E\in\mathcal{T}_h} \Big\Vert(\Pi_{1}^{0,E}\frac{\partial (\mathbf{u}_h)_i}{\partial x_j}-\frac{\partial (\mathbf{u}_h)_i}{\partial x_j}\Big\Vert_{0,E}^2  \Vert(\Pi_2^{0,E}\mathbf{u}_h)_j \Vert_{L^\infty(E)}^2\Big)^{\frac{1}{2}} \vert\mathbf{\Phi}\vert_{1,\Omega}\\
     & \lesssim 
    \Big(\sum_{E\in\mathcal{T}_h}\|\mathbf{\Pi}_{k-1}^{0,E}\nabla\mathbf{u}_h-\nabla\mathbf{u}_h\|_{0,E}^2\|\mathbf{\Pi}_k^{0,E}\mathbf{u}_h\|^2_{L^\infty(E)}\Big)^{1/2}\vert\mathbf{\Phi}\vert_{1,\Omega} \,.
    \end{split}
    \end{equation*}
The term $\|\mathbf{\Pi}_{k-1}^{0,E}\nabla\mathbf{u}_h-\nabla\mathbf{u}_h\|_{0,E}$ can be bounded as follows:
\begin{equation*}
\begin{split}
     \|\mathbf{\Pi}_{k-1}^{0,E}\nabla\mathbf{u}_h-\nabla\mathbf{u}_h\|_{0,E} &\leq \Vert\mathbf{\Pi}_{k-1}^{0,E}\nabla\mathbf{u}_h-\nabla \mathbf{\Pi}^{\nabla,E}_k\mathbf{u}_h \Vert_{0,E}+
    \Vert\nabla (\mathbf{\Pi}^{\nabla,E}_k\mathbf{u}_h-\mathbf{u}_h)\Vert_{0,E}\\
    & \hskip -0.7cm =\Vert\mathbf{\Pi}_{k-1}^{0,E}\nabla\mathbf{u}_h-\mathbf{\Pi}_{k-1}^{0,E}\nabla \mathbf{\Pi}^{\nabla,E}_k\mathbf{u}_h \Vert_{0,E}
   +  \Vert\nabla (\mathbf{\Pi}^{\nabla,E}_k\mathbf{u}_h-\mathbf{u}_h)\Vert_{0,E}\\
    & \hskip -0.7cm  \lesssim 2
    \Vert\nabla (\mathbf{\Pi}^{\nabla,E}_k\mathbf{u}_h-\mathbf{u}_h) \Vert_{0,E}\lesssim \sigma_E \,.
\end{split}
\end{equation*}
We conclude that
\begin{equation}\label{eq:bound-T2}
    \begin{split}
    (T_2)& \lesssim \Big( \sum_{E\in\mathcal{T}_h} \sigma_E^2\|\mathbf{\Pi}_k^{0,E}\mathbf{u}_h\|^2_{L^\infty(E)}\Big)^{1/2}\vert\mathbf{\Phi}\vert_{1,\Omega} = \Big(\sum_{E\in\mathcal{T}_h}\eta^2_{c2,E}\Big)^{1/2}\vert\mathbf{\Phi}\vert_{1,\Omega} \,.
\end{split}
\end{equation}

Next, considering $(T_3)$, we start again by using H\"older's inequality with exponents $\left(\frac{1}{2},\frac{1}{4},\frac{1}{4}\right)$ to get
\begin{equation*} 
      (T_3)\lesssim  \sum_{i,j=1}^2\sum_{E\in\mathcal{T}_h} \Big\Vert\frac{\partial (\mathbf{u}_h)_i}{\partial  x_j}\Big\Vert_{0,E}   \|(\mathbf{\Pi}_k^{0,E}\mathbf{u}_h-\mathbf{u}_h)_j\|_{L^4(E)}\|\mathbf{\Phi}_i\|_{L^4(E)}  \,.
\end{equation*}
Now,  we apply to the quantity $\psi = (\mathbf{\Pi}_k^{0,E}\mathbf{u}_h-\mathbf{u}_h)_j$ the inverse inequality $\Vert \psi \Vert_{L^4(E)} \lesssim h_E^{-1/2} \Vert \psi \Vert_{0,E}$ and the scaled Poincar\'e inequality for zero-mean functions $\Vert \psi \Vert_{0,E} \lesssim h_E \Vert \nabla \psi \Vert_{1,E}$, whence $\Vert \psi \Vert_{L^4(E)} \lesssim h_E^{1/2} \Vert \nabla \psi \Vert_{1,E}$. Invoking as for $(T_2)$ the Sobolev embedding $ H^1(E) \subset L^4(E)$ for $\boldsymbol{\Phi}_i$, we obtain
\begin{equation*}
     (T_3)  \lesssim   \Big(\sum_{E\in\mathcal{T}_h} \| \nabla\mathbf{u}_h\|_{0,E}^2\|\nabla(\mathbf{\Pi}_k^{0,E}\mathbf{u}_h-\mathbf{u}_h) \|^2_{0,E}\Big)^{1/2}\vert\mathbf{\Phi}\vert_{1,\Omega} \,.
\end{equation*}
Writing
$$
\| \nabla\mathbf{u}_h\|_{0,E} = \| \nabla\mathbf{u}_h  - \nabla \mathbf{\Pi}_k^{\nabla,E}\mathbf{u}_h +  \nabla \mathbf{\Pi}_k^{\nabla,E}\mathbf{u}_h   \|_{0,E} \lesssim \sigma_E + \Vert \nabla \mathbf{\Pi}_k^{\nabla,E}\mathbf{u}_h \Vert_{0,E} 
$$
and
\begin{equation*}
\begin{split}
\|\nabla(\mathbf{\Pi}_k^{0,E}\mathbf{u}_h-\mathbf{u}_h) \|_{0,E} & =  \|\nabla(\mathbf{\Pi}_k^{0,E}\mathbf{u}_h-
\mathbf{\Pi}_k^{\nabla,E}\mathbf{u}_h +  \mathbf{\Pi}_k^{\nabla,E}\mathbf{u}_h -  \mathbf{u}_h) \|_{0,E}  \\
& \lesssim \sigma_E +  \|\nabla(\mathbf{\Pi}_k^{0,E}\mathbf{u}_h  - \mathbf{\Pi}_k^{\nabla,E}\mathbf{u}_h ) \|_{0,E} \,,
      \end{split}
\end{equation*}
we conclude that
\begin{equation}\label{eq:bound-T3}
    (T_3) \lesssim \Big(\sum_{E\in\mathcal{T}_h}\eta^2_{c3,E}\Big)^{1/2}\vert\mathbf{\Phi}\vert_{1,\Omega} \,.
\end{equation}

At last, we analyze term $J_7$ which can be written as
\begin{equation}
    J_7=\sum_{E\in\mathcal{T}_h}\Big(\mathbf{t}^E_h,\mathbf{\Phi}-\mathbf{\Phi}_I-\mathbf{\Pi}_k^{0,E}(\mathbf{\Phi}-\mathbf{\Phi}_I)\Big)_E \,.
\end{equation}
Using similar calculations as done for $J_6$ (see \eqref{eq:split_j7}),  we write 
\begin{equation}
\begin{split}
        J_7&=\sum_{E\in\mathcal{T}_h}\Big(\mathbf{t}^E_h-\mathbf{\Pi}_k^{0,E}\mathbf{t}^E_h,\mathbf{\Phi}-\mathbf{\Phi}_I\Big)_E\\
        &\leq\sum_{E\in\mathcal{T}_h} \Vert\mathbf{t}^E_h-\mathbf{\Pi}_k^{0,E}\mathbf{t}^E_h\Vert_{0,E}\left\Vert\mathbf{\Phi}-\mathbf{\Phi}_I\right\Vert_{0,E}
        \lesssim\sum_{E\in\mathcal{T}_h} \Vert\mathbf{t}^E_h-\mathbf{\Pi}_k^{0,E}\mathbf{t}^E_h \Vert_{0,E}h_E\vert\mathbf{\Phi}\vert_{1,D_E}\\
        &\lesssim\Big(\sum_{E\in\mathcal{T}_h}h_E^2 \Vert\mathbf{t}^E_h-\mathbf{\Pi}_k^{0,E}\mathbf{t}^E_h\Vert^2_{0,E}\Big)^{1/2}\Big(\sum_{E\in\mathcal{T}_h}\vert\mathbf{\Phi}\vert_{1,D_E}^2\Big)^{1/2}\\
        &\lesssim \Big(\sum_{E\in\mathcal{T}_h}h_E^2 \Vert\mathbf{t}^E_h-\mathbf{\Pi}_k^{0,E}\mathbf{t}^E_h\Vert^2_{0,E}\Big)^{1/2} \left\vert\mathbf{\Phi}\right\vert_{1,\Omega}  \,.
\end{split}
\end{equation}
So we obtain the same upper bound as for the term $(T_1)$, namely
\begin{equation}\label{eq:Bound-J7}
    J_7\lesssim\Big(\sum_{E\in\mathcal{T}_h}\eta^2_{c1,E}\Big)^{1/2}\vert\mathbf{\Phi}\vert_{1,\Omega} \,. \end{equation}

Collecting the bounds \eqref{eq:BD-J1}-\eqref{eq:BD-J3}, \eqref{eq:bound-J6}-\eqref{eq:bound-T3} and \eqref{eq:Bound-J7}, and recalling the identity \eqref{eq:phi-identity}, we 
obtain
$$
a(\mathbf{\Phi},\mathbf{\Phi}) = \nu \vert \mathbf{\Phi} \vert_{1,\Omega}^2 \lesssim \eta \vert \mathbf{\Phi} \vert_{1,\Omega}\,
$$
whence $\nu \vert \mathbf{\Phi} \vert_{1,\Omega} \lesssim \eta$, i.e. $a(\mathbf{\Phi},\mathbf{\Phi}) \lesssim \frac1\nu \eta^2$.
We conclude our proof, thanks to the bound \eqref{eq:errbnd-1}. 
\end{proof}

\subsection{Lower error bounds}\label{sec:lower-bounds}

While a standard local a posteriori lower bound of the error, similar e.g. to \cite[Theorem 16]{a_post_cangiani} for VEM discretizations, does not seem to be easily provable, we are able to establish partial results, collected in the two following propositions.

The first proposition states that the bulk and jump estimators provide a lower bound for a local energy error (scaled differently with respect to $\nu$ if compared to \eqref{eq:est_post}), up to data oscillation and a stabilization term.

\begin{proposition}[Lower bounds involving the bulk and jump estimators]\label{prop:lower-bounds-1} For any $E \in {\cal T}$, it holds
\begin{eqnarray}
\eta_{B,E} &\lesssim&  \nu\Vert\mathbf{u}-\mathbf{u}_h\Vert_{1,E}+\|p-p_h\|_{0,E} + \eta_{f,E} + (\nu^{1/2}+\nu)\sigma_E \,; \label{eq:AAA}\\
\eta_{e,E} &\lesssim&   \nu\Vert\mathbf{u}-\mathbf{u}_h\Vert_{1,D_E}+\|p-p_h\|_{0,D_E} + \!\! \sum_{E' \in D_E} \!\!\!\! \big( \eta_{f,E'} + (\nu^{1/2}+\nu)\sigma_{E'} \big). \label{eq:BBB}
\end{eqnarray}
 \end{proposition}
 \begin{proof} The inequalities can be obtained by using Verf\"urth's classical technique based on localized bubble functions (see \cite{verfurth96}; see also \cite{kanschat2008,a_post_cangiani}). 
 We refrain from reporting all the details of the lengthy derivation, and we just highlight the most peculiar elements of the analysis.
 
 In order to prove \eqref{eq:AAA}, let us introduce the piecewise polynomial bubble function $b_E \in H^1_0(E)$ which satisfies for any polynomial $\varphi \in \mathbb{P}_{2k}(E)$ the bounds
 \begin{equation}\label{eq:prop-bubble}
 \begin{split}
  &\Vert \varphi \Vert_{0,E} \lesssim \Vert \varphi b_E^{1/2} \Vert_{0,E} \,, \qquad  \Vert \varphi b_E \Vert_{0,E} \lesssim \Vert \varphi  \Vert_{0,E} \,,  \\
  & \Vert \nabla (\varphi b_E) \Vert_{0,E} \lesssim h_E^{-1} \Vert \varphi  \Vert_{0,E} \,, \qquad  \Vert \varphi b_E \Vert_{L^\infty(E)} \lesssim  h_E^{-1}\Vert \varphi  \Vert_{0,E} \,. 
 \end{split}
 \end{equation}
 Setting $\mathbf{R}_E =\mathbf{f}_h+\nu\Delta\mathbf{\Pi}^{\nabla,E}_k\mathbf{u}_h+\nabla p_h- (\mathbf{\Pi}^{0,E}_{k-1}\nabla\mathbf{u}_h) \mathbf{\Pi}^{0,E}_k\mathbf{u}_h$  
 and $\boldsymbol{\varrho}_E=\mathbf{R}_E b_E$, one has
 \begin{equation*}
 \eta_{B,E}^2 = h_E^2 \Vert \mathbf{R}_E \Vert_{0,E}^2 \lesssim h_E^2 (\mathbf{R}_E, \boldsymbol{\varrho}_E)_{0,E} 
 \end{equation*} 
 with
 \begin{equation}\label{eq:resR}
 \begin{split}
 (\mathbf{R}_E, \boldsymbol{\varrho}_E)_{0,E} &= (\mathbf{f}_h - \mathbf{f}, \boldsymbol{\varrho}_E)_{0,E} + \nu (\nabla \mathbf{u} - \nabla \mathbf{\Pi}^{\nabla,E}_k \mathbf{u}_h, \nabla \boldsymbol{\varrho}_E)_{0,E} +(p-p_h, \nabla \cdot \boldsymbol{\varrho}_E)_{0,E} \\
 & \qquad - ( \underbrace{(\nabla\mathbf{u}) \mathbf{u} -  (\mathbf{\Pi}^{0,E}_{k-1}\nabla\mathbf{u}_h) \mathbf{\Pi}^{0,E}_k\mathbf{u}_h}_{\mathbf{NL}} , \boldsymbol{\varrho}_E)_{0,E}\,.
 \end{split}
 \end{equation}
 Focussing just on the nonlinear term, it can be split telescopically as
\begin{equation*}
 \begin{split}
\mathbf{NL} &= \Big( (\nabla\mathbf{u}) \mathbf{u} - (\nabla\mathbf{u}) \mathbf{u}_h \Big)  +  \Big( (\nabla\mathbf{u}) \mathbf{u}_h - (\nabla\mathbf{u}_h) \mathbf{u}_h \Big) \\
& \ + \Big( (\nabla\mathbf{u}_h) \mathbf{u}_h - (\nabla\mathbf{u}_h) \mathbf{\Pi}^{0,E}_k\mathbf{u}_h\Big) + \Big( (\nabla\mathbf{u}_h) \mathbf{\Pi}^{0,E}_k\mathbf{u}_h -
(\mathbf{\Pi}^{0,E}_{k-1}\nabla\mathbf{u}_h) \mathbf{\Pi}^{0,E}_k\mathbf{u}_h \Big)\,.
 \end{split}
 \end{equation*}
 Considering the first addend, we have
 \begin{equation*}
 \begin{split}
 \vert ((\nabla\mathbf{u}) \mathbf{u} - (\nabla\mathbf{u}) \mathbf{u}_h, \boldsymbol{\varrho}_E)_{0,E} \vert  & \lesssim 
 \Vert \nabla\mathbf{u} \Vert_{0,E} \Vert \mathbf{u} - \mathbf{u}_h \Vert_{0,E} \Vert \boldsymbol{\varrho}_E \Vert_{L^\infty(E)} \\
 & \lesssim \nu^{-1} \Vert \mathbf{f} \Vert_{-1,\Omega} \Vert \mathbf{u} - \mathbf{u}_h \Vert_{1,E} h_E^{-1} \Vert \mathbf{R}_E \Vert_{0,E} \\
 & \lesssim \nu  \Vert \mathbf{u} - \mathbf{u}_h \Vert_{1,E} h_E^{-1} \Vert \mathbf{R}_E \Vert_{0,E} \,,
 \end{split}
 \end{equation*}
 where we have used \eqref{eq:impl_small_data}, \eqref{gamma_small} and \eqref{eq:prop-bubble}. The remaining addends of $\mathbf{NL}$ can be bounded similarly, using now \eqref{eq:u_err} and \eqref{eq:as_gamma_h_NS}, and involving the stabilization term $\sigma_E$ in estimating the projection errors, as done in the proof of Theorem \ref{theo:reliability}.
 
 The proof of \eqref{eq:BBB} follows a similar procedure: one uses a bubble supported about an edge $e$, and by an integration by parts one is led to bound integrals similar to the previous ones over the elements sharing the edge. 
 \end{proof}
 
 Next results concern the estimators $\eta_{c1,E}, \, \eta_{c2,E}, \, \eta_{c3,E}$ that were introduced to control the error in the nonlinear convective term.  Note that $\eta_{c1,E}$ measures a projection error for the quantity $(\mathbf{\Pi}_{k-1}^{0,E}\nabla\mathbf{u}_h) \mathbf{\Pi}_k^{0,E}\mathbf{u}_h$; hence, it cannot be bounded from above in terms of the velocity discretization error and the stabilization term only (which could both vanish without the vanishing of the projection error). For this reason, we also make the projection error for  $(\nabla\mathbf{u})\mathbf{u}$ appear in the estimate.
 
 \begin{proposition}[Lower bounds involving the convective-term estimators]\label{prop:lower-bounds-1} For any $E \in {\cal T}$, it holds
\begin{eqnarray}
\eta_{c1,E} &\lesssim& h_E \Vert (\nabla\mathbf{u})\mathbf{u} - \mathbf{\Pi}_{k}^{0,E} \big( (\nabla\mathbf{u}) \mathbf{u}  \big) \Vert_{0,E} + \nu\Vert\mathbf{u}-\mathbf{u}_h\Vert_{1,E} +\nu \sigma_E
 \,;  \label{eq:bound-etac1} \\ 
\eta_{c2,E} &\lesssim&  (1+ \vert \log h_E \vert) \nu \sigma_E \,; \label{eq:bound-etac2}\\
\eta_{c3,E} &\lesssim&  (\nu + \sigma_E) \sigma_E \label{eq:bound-etac3}\,. 
\end{eqnarray}
 \end{proposition}
 \begin{proof} 
 We start by establishing \eqref{eq:bound-etac1}. To this end, it is convenient to set $\mathbf{d}=(\nabla\mathbf{u})\mathbf{u}$ and $\mathbf{d}_h = (\mathbf{\Pi}_{k-1}^{0,E}\nabla\mathbf{u}_h) \mathbf{\Pi}_k^{0,E}\mathbf{u}_h$, so that $\eta_{c1,E}=h_E \Vert \mathbf{d}_h - \mathbf{\Pi}_k^{0,E} \mathbf{d}_h \Vert_{0,E}$. Since $\mathbf{d}_h \in (\mathbb{P}_{k(k-1)}(E))^2$, we can write
 $$
 \Vert \mathbf{d}_h - \mathbf{\Pi}_k^{0,E} \mathbf{d}_h \Vert_{0,E} = \sup_{\boldsymbol{\varphi} \in (\mathbb{P}_{k(k-1)}(E))^2} \frac{(\mathbf{d}_h - \mathbf{\Pi}_k^{0,E} \mathbf{d}_h, \boldsymbol{\varphi})_{0,E}}{\Vert \boldsymbol{\varphi} \Vert_{0,E}}\,,
 $$
 and the numerator can be written as
 \begin{equation*}
 \begin{split}
 (\mathbf{d}_h - \mathbf{\Pi}_k^{0,E} \mathbf{d}_h, \boldsymbol{\varphi})_{0,E} &= (\mathbf{d} - \mathbf{\Pi}_k^{0,E} \mathbf{d}, \boldsymbol{\varphi})_{0,E} \\
 & \ \ + (\mathbf{d}_h - \mathbf{d}, \boldsymbol{\varphi})_{0,E} - (\mathbf{d}_h - \mathbf{d}, \mathbf{\Pi}_k^{0,E} \boldsymbol{\varphi})_{0,E}\,.
 \end{split}
 \end{equation*}
 The second term can be bounded as
 $$
 \vert (\mathbf{d}_h - \mathbf{d}, \boldsymbol{\varphi})_{0,E} \vert \leq \Vert \mathbf{d} - \mathbf{d}_h  \Vert_{L^1(E)} \Vert \boldsymbol{\varphi} \Vert_{L^\infty(E)} \lesssim h_E^{-1} \Vert \mathbf{d} - \mathbf{d}_h \Vert_{L^1(E)} \Vert \boldsymbol{\varphi} \Vert_{0,E}\,, 
 $$
 and similarly for the third term. Thus,
 $$
 \Vert \mathbf{d}_h - \mathbf{\Pi}_k^{0,E} \mathbf{d}_h \Vert_{0,E} \lesssim   \Vert \mathbf{d} - \mathbf{\Pi}_k^{0,E} \mathbf{d}
 \Vert_{0,E}  + h_E^{-1} \Vert \mathbf{d} - \mathbf{d}_h \Vert_{L^1(E)}\,.
 $$
 Finally, we observe that  $\mathbf{d} - \mathbf{d}_h$ is precisely the term $\mathbf{NL}$ introduced in \eqref{eq:resR}, which has been already bounded above, leading to \eqref{eq:bound-etac1}.
 
 In order to prove \eqref{eq:bound-etac2}, we invoke inverse inequalities to get
 \begin{equation*}
 \begin{split}
 \Vert\mathbf{\Pi}_k^{0,E}\mathbf{u}_h \Vert_{L^\infty(E)} & \leq  \Vert \mathbf{u}_h - \mathbf{\Pi}_k^{0,E}\mathbf{u}_h \Vert_{L^\infty(E)} +  \Vert \mathbf{u}_h  \Vert_{L^\infty(E)} \\
 & \lesssim h_E^{-1} \Vert \mathbf{u}_h - \mathbf{\Pi}_k^{0,E}\mathbf{u}_h \Vert_{0,E} +  h_E^{-2/q} \Vert \mathbf{u}_h  \Vert_{L^q(E)} \qquad \forall q \in (2, +\infty)\,, \\ 
 & \lesssim \vert \mathbf{u}_h \vert_{1,E} + h_E^{-2/q} \Vert \mathbf{u}_h  \Vert_{L^q(\Omega)} \\
  & \lesssim \vert \mathbf{u}_h \vert_{1,\Omega} + \frac{q}2 h_E^{-2/q} \Vert \mathbf{u}_h  \Vert_{1,\Omega} 
  \leq (1 + \frac{q}2  \, {\rm e}^{\frac2{q} \vert \log h_E \vert}) \Vert \mathbf{u}_h  \Vert_{1,\Omega}  \,, 
 \end{split} 
 \end{equation*}
 where all the constants implied by the symbol $\lesssim$ can be chosen independent of $q$. Thus, we can take $q = 2  \vert \log h_E \vert$, which gives the desired result.
 
 Finally, the bound \eqref{eq:bound-etac3} immediately follows from the by now familiar bounds $\Vert\nabla(\mathbf{\Pi}_k^{0,E}\mathbf{u}_h-\mathbf{\Pi}_k^{\nabla,E}\mathbf{u}_h)\Vert_{0,E} \lesssim \sigma_E$ and $\Vert\nabla\mathbf{\Pi}^{\nabla,E}_k\mathbf{u}_h \Vert_{0,E}  \lesssim \nu$.
 \end{proof}

%%%%%%%%%%%%%%%%%%%%%%%%%%%%%%%%%
\section{Numerical results}
\label{sec:experiments}

Throughout this section, we choose the order of accuracy $k=2$ in the definition of the VEM spaces. 

\begin{figure} [!htb]
    \centering
    \includegraphics[scale=0.51]{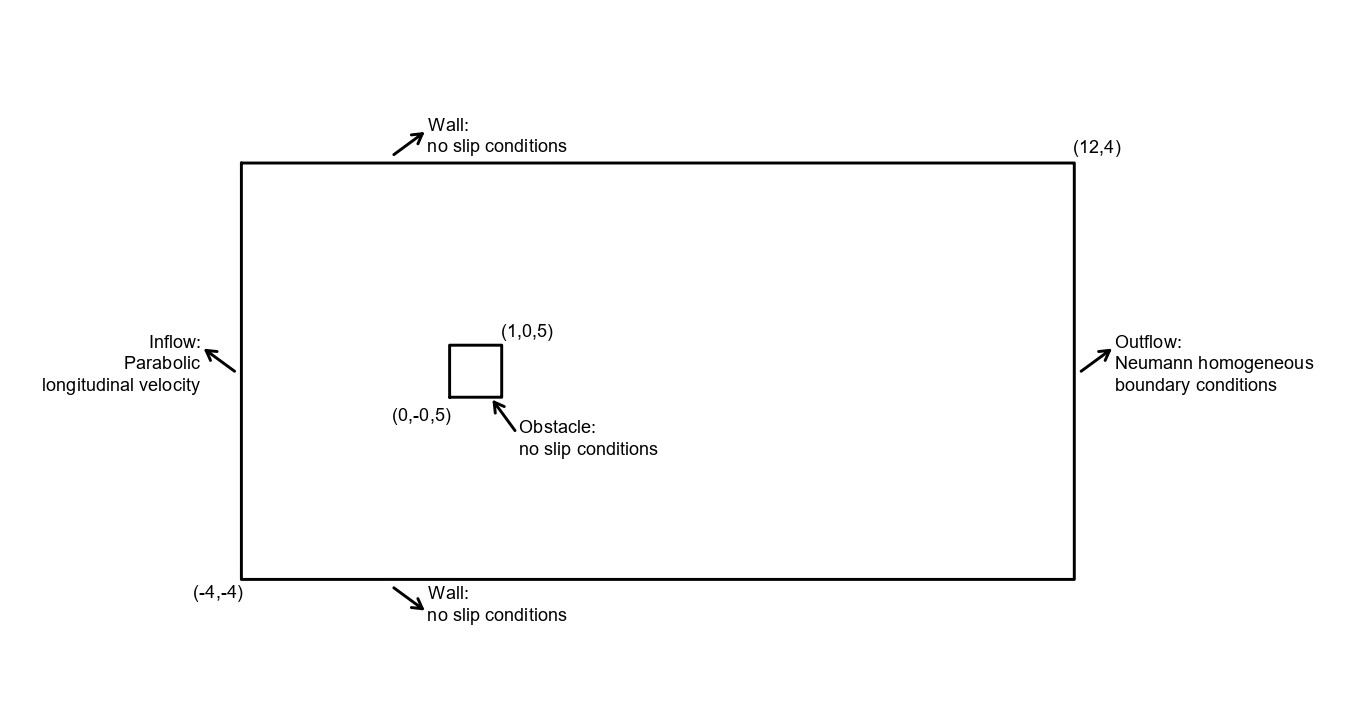}
    \caption{Domain of interest with the boundary conditions of the fluid-dynamics test case}
    \label{fig:bound_fluid_imm}
\end{figure}

To assess the performance of our error estimator, we consider the domain shown in Figure \ref{fig:bound_fluid_imm}, namely a rectangular channel with a square obstacle whose centroid is placed on the longitudinal axis; this can be viewed as the horizontal section of a 3D channel containing an infinite square cylinder in the middle, and has important applications, e.g., in Wind Engineering (see  \cite{Fluid_1}). 
We will use both manufactured (i.e., analytical) solutions to monitor the decay of the true error and the estimator, and the (unknown) solution of a fluid-dynamics problem to investigate the effect of mesh-refinement driven by our estimator.

All our meshes are obtained by uniform or adaptive refinements from an initial, uniform partition of the domain into squares with edges equal to half the edge of the square cylinder. The classical quad-tree refinement is applied to each refined square, which replaces it by four equal sub-squares; this may generate hanging nodes. Following the VEM philosophy, a square containing hanging nodes is viewed as a polygon. For easiness of implementation, we confine ourselves to refinements with a maximum of one hanging node per edge (although from the theoretical point of view we could choose any bound for the maximum number of hanging nodes per edge). Thus, our VEM elements may range from pentagons (squares with one hanging node) to octagons (squares with four hanging nodes).

In some of the forthcoming plots, we will monitor the different components of the a posteriori estimator, which will be denoted by $\eta_s$ with
\begin{equation}
    \eta_s^2=\sum_{E\in\mathcal{T}_h}\eta_{s,E}^2 \,,
\end{equation}
where the subscript $s$ stands for one of the
subscripts previously used in the definition \eqref{eq:def-components}.
We point out that we will enforce non-homogeneous Dirichlet boundary conditions and 
homogeneous Neumann boundary conditions, namely 
\begin{equation}
    \begin{cases}
    \mathbf{u}=\mathbf{g}\quad \text{on}\quad\partial\Omega_D\,, \\
    (\nu\mathbf{\nabla}\mathbf{u}+p\mathbf{I})\mathbf{n}=\mathbf{0}
    \quad \text{on} \quad \partial\Omega_N \,,\\
    \end{cases}
\end{equation}
even if in Sect. \ref{sec:estimate} the estimator was defined assuming homogeneous Dirichlet boundary conditions only. The extension to the general case poses no difficulty.

The code used for implementation was written in \textit{Matlab}, taking inspiration from \cite{Sutton} and \cite{Guide}.

In the following subsections we will briefly describe Newton's method to solve the nonlinear problem, and the adaptive refinement algorithm.  The last part of the section is devoted to the computational results.

\subsection{Newton's method}
To solve the discrete Navier-Stokes equations \eqref{eq:NS_discr}, we used Newton's method as described in \cite{Newton}. Starting from an initial guess $\mathbf{u}_0\in\mathbf{V}_h$, the Newton iterates $\{\mathbf{u}_m\in\mathbf{V}_h,p_m\in Q_h\}_{m=1,2,...}$ are defined by
\begin{equation}\label{newton_iter}
\begin{split}
        a_h(\mathbf{u}_m,\mathbf{v}_h)+{c}_h(\mathbf{u}_m,\mathbf{u}_{m-1},\mathbf{v}_h)+{c}_h(\mathbf{u}_{m-1},\mathbf{u}_{m},\mathbf{v}_h)+b(\mathbf{v}_h,p_m)+b(\mathbf{u}_m,q)\\
        =(\mathbf{f}_h,\mathbf{v}_h)+{c}_h(\mathbf{u}_{m-1};\mathbf{u}_{m-1},\mathbf{v}_h)\quad\forall\, \mathbf{v}_h\in\mathbf{V}_h,q\in Q_h \,.
\end{split}
\end{equation}
The initial guess $\mathbf{u}_0$ is chosen as the solution of the corresponding Stokes problem, obtained as described in  \cite{Stokes_1}.
As a stopping rule,  we impose the condition $\| \bar{\mathbf{u}}_m-\bar{\mathbf{u}}_{m-1}\|_{\ell_2 }<\textit{tol}$, where $\bar{\mathbf{u}}_m$ is the array containing the values of the degrees of freedom of ${\mathbf{u}}_m$. The tolerance $\textit{tol}$ is a small arbitrary value chosen equal to $10^{-9}$. The maximum number of iterates is taken equal to $10$.

\subsection{Refinement algorithm}

The algorithm used to refine the mesh $\mathcal{T}_h$ was the classical $\texttt{SOLVE} \rightarrow
    \texttt{ESTIMATE}
    \rightarrow
    \texttt{MARK}
    \rightarrow
    \texttt{REFINE}$ loop. 
    The subset of marked elements ${\cal M}_h$ is chosen according to D\"orfler's recipe
\begin{equation}
    \sum_{E\in{\mathcal{M}_h}}\eta_E^2\geq\theta \sum_{E\in {\mathcal{T}_h}}\eta_E^2 \,,
\end{equation}
for some fixed $\theta \in (0,1)$. 
It is important to underline that, generally, the previous rule might produce meshes with an arbitrary number of hanging nodes per edge. Therefore, we implemented an additional recursive refinement procedure to satisfy the constraint of having at most one hanging node per edge.

\subsection{Numerical tests with analytical solutions}
\begin{figure}[t!]
        \centering
        \begin{subfigure}[b]{0.5\textwidth}
                \centering
         \includegraphics[width=\textwidth]{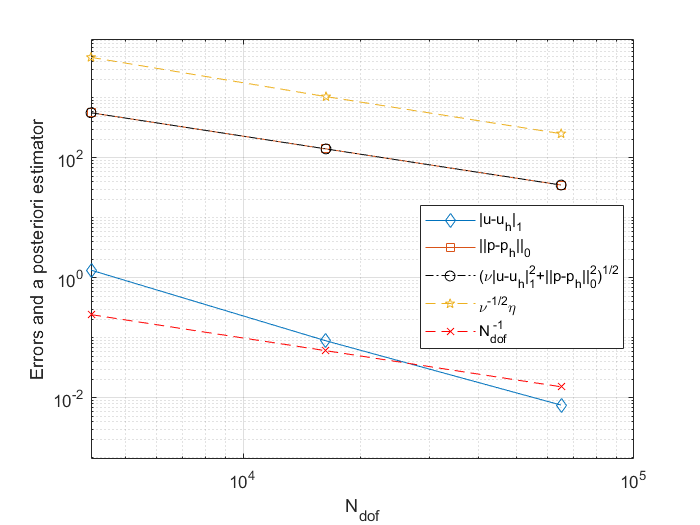}
         \caption{$Re$=1}
        \end{subfigure}%    <-- % added here
        \hfill %% useful if width of each figure is less the .5\textwidth
        \begin{subfigure}[b]{0.5\textwidth}
                \centering
         \includegraphics[width=\textwidth]{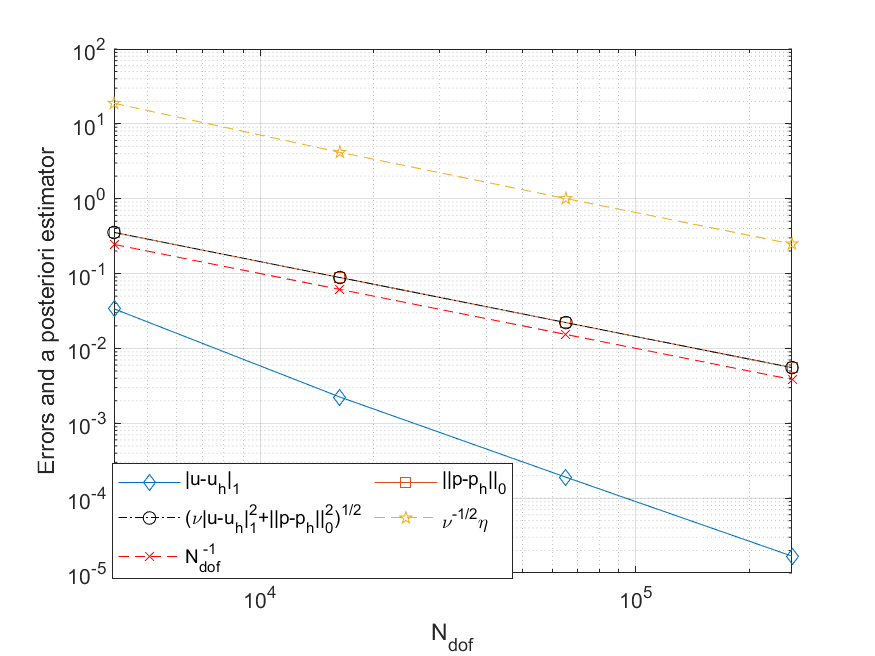}
         \caption{$Re$=40}
        \end{subfigure}
        \caption{Test 1. Uniform refinements: errors and a posteriori estimator}
        \label{fig:unif_1_err}
\end{figure}
\begin{figure}[t!]
        \centering
        \begin{subfigure}[b]{0.5\textwidth}
                \centering
         \includegraphics[width=\textwidth]{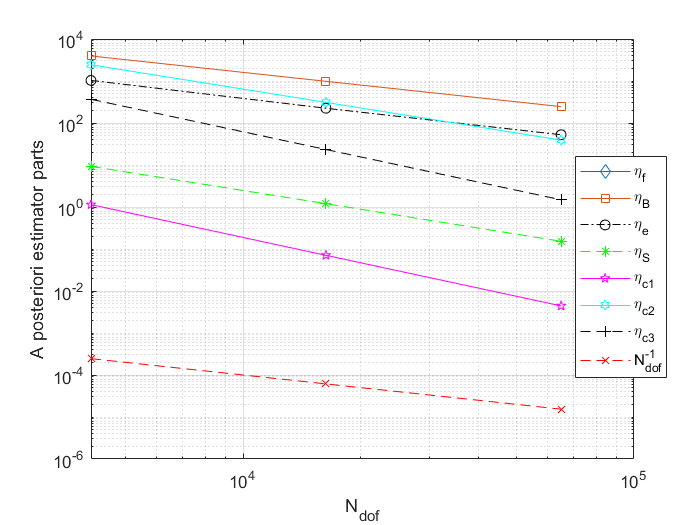}
         \caption{$Re$=1}
        \end{subfigure}%    <-- % added here
        \hfill %% useful if width of each figure is less the .5\textwidth
        \begin{subfigure}[b]{0.5\textwidth}
                \centering
         \includegraphics[width=\textwidth]{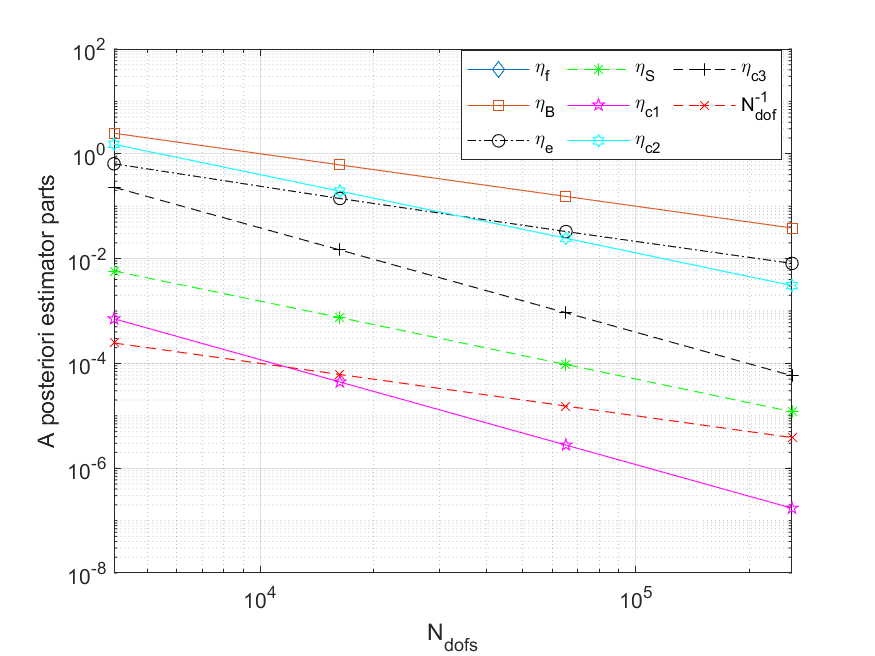}
         \caption{$Re$=40}
        \end{subfigure}
        \caption{Test 1. Components of the a posteriori estimator during uniform refinements}
        \label{fig:unif_1_parts}
\end{figure}
Hereafter we report the results of some numerical experiments with manufactured solutions.
We underline that in the plots related to velocity errors,  we indicate such errors as $\vert\mathbf{u}-\mathbf{u}_h\vert_1$ but we actually mean $\left(\sum_{E\in\mathcal{T}_h}\|\mathbf{\nabla}\mathbf{u}-\mathbf{\Pi}^{0,E}_1\mathbf{\nabla}\mathbf{u}_h\|^2_{0,E}\right)^{1/2}$, since the explicit expression of the VEM solution is not known.

\subsubsection{Test 1}

Taking inspiration from \cite{Navier-Stokes_1}, we choose the following solution:
\begin{equation}
    \mathbf{u}=\nu\begin{bmatrix}
3(x^2-y^2) \\
-2xy
\end{bmatrix}\,, \qquad
p=\nu^2(\frac{9}{2}(x^2+y^2)^2-C) \,,
\end{equation}
where $C=\int_\Omega\frac{9}{2}(x^2+y^2)^2\,d\Omega/\vert\Omega\vert$. Dirichlet boundary conditions equal to the values of exact solution are enforced on $\partial\Omega$. With this solution one has $\mathbf{f}=\mathbf{0}$, which implies $\eta_f=0$.  We applied uniform refinements of the initial mesh for different values of $Re$. In Figure \ref{fig:unif_1_err} we show, in log-log scale, the behavior of the errors and the estimator for $Re=1$ and $Re=40$ vs the number of active degrees of freedom, whereas in Figure \ref{fig:unif_1_parts} we show the different components of the a posteriori estimator. We notice that the bulk error $\eta_B$ is the dominant one. We underline that the estimator and the pressure error follow the optimal convergence rate (indicated by the dashed line with symbol $\text{N}_\text{dof}^{-1}$), whereas the velocity error seems to decrease faster while adding degrees of freedoms. Note that the values of the estimator and the error in the case $Re=1$ are higher than the corresponding values for $Re=40$, simply because the exact solution is larger in the former case than in the latter. In this test Newton's method converges in 3 iterations or less.

\subsubsection{Test 2}
We now consider the $\nu$-independent solution
\begin{equation}
    \mathbf{u}=\begin{bmatrix}
e^{-(x-12)/6}\sin(y/6)\\
-e^{-(x-12)/6}\cos(y/6)
\end{bmatrix}\,, \qquad
p=e^{(-(x-12)/6)}\sin(y/6) \,,
\end{equation}
which satisfies $\nu\boldsymbol{\Delta}u=\mathbf{0}$, whence $\mathbf{f}=(\boldsymbol{\nabla}\mathbf{u})\mathbf{u}-\nabla p$. 
As in \textit{Test 1}, we enforce Dirichlet conditions on $\partial\Omega$ and we use uniform mesh refinements.
\begin{figure}[t]
        \centering
        \begin{subfigure}[t]{0.5\textwidth}
                \centering
         \includegraphics[width=\textwidth]{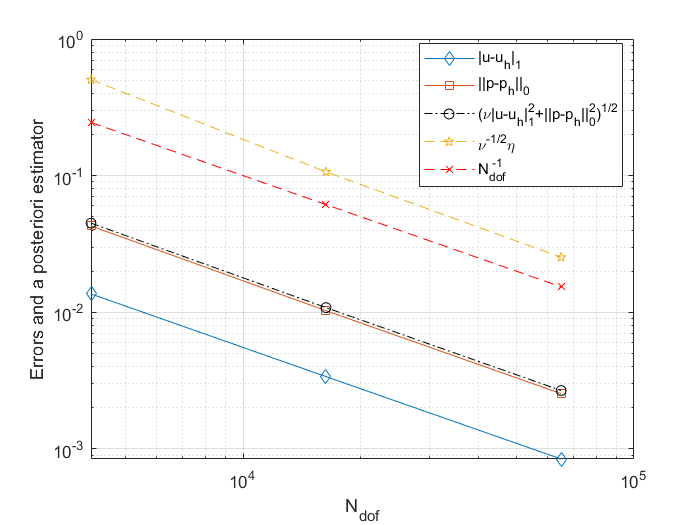}
         \caption{$Re$=1}
        \end{subfigure}%    <-- % added here
        \hfill %% useful if width of each figure is less the .5\textwidth
        \begin{subfigure}[t]{0.5\textwidth}
                \centering
         \includegraphics[width=\textwidth]{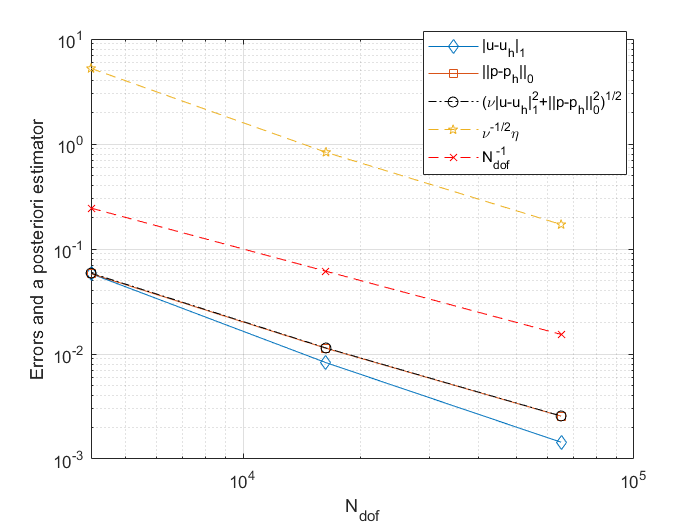}
         \caption{$Re$=40}
        \end{subfigure}
        \caption{Test 2. Errors and a posteriori estimator during uniform refinements}
         \label{fig:unif_2_err}
\end{figure}
\begin{figure}[t]
        \centering
        \begin{subfigure}[t]{0.5\textwidth}
                \centering
         \includegraphics[width=\textwidth]{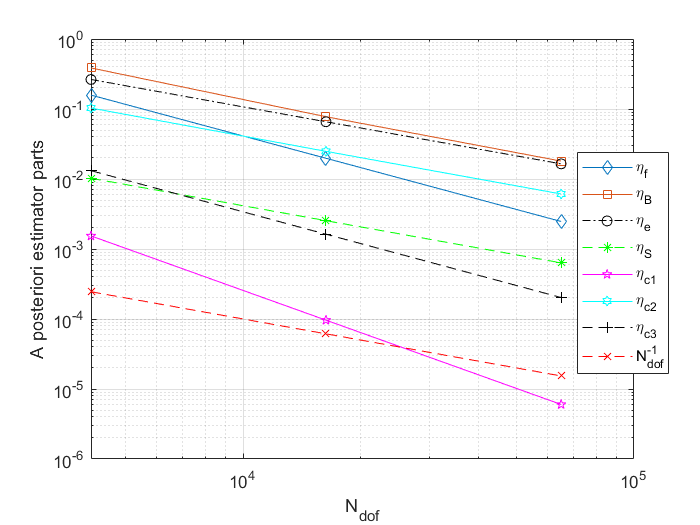}
         \caption{$Re$=1}
        \end{subfigure}%    <-- % added here
        \hfill %% useful if width of each figure is less the .5\textwidth
        \begin{subfigure}[t]{0.5\textwidth}
                \centering
         \includegraphics[width=\textwidth]{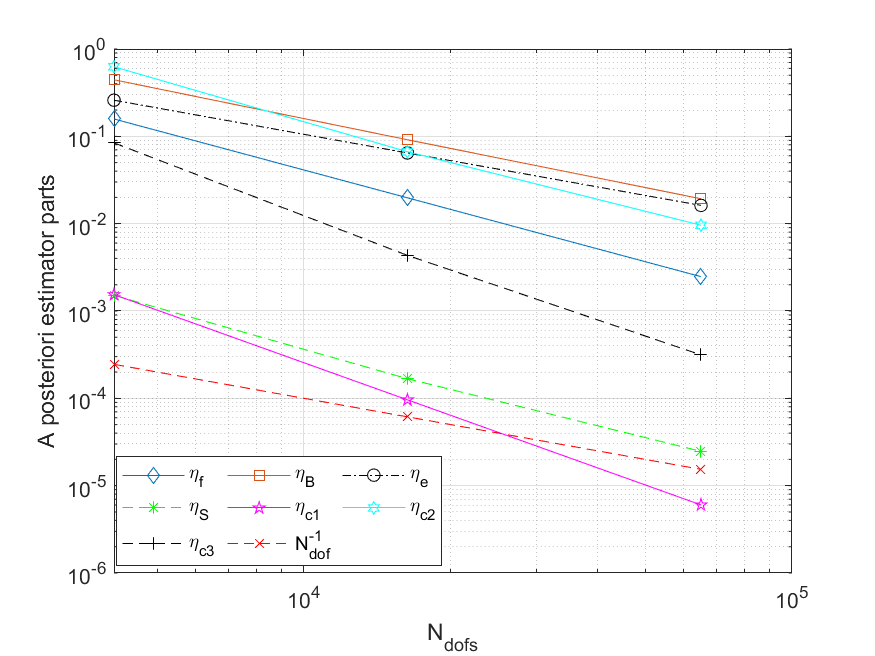}
         \caption{$Re$=40}
        \end{subfigure}
         \caption{Test 2. Components of the a posteriori estimator during uniform refinements}
         \label{fig:unif_2_parts}
\end{figure}
 In Figure \ref{fig:unif_2_err} we report the errors and the a posteriori estimator vs the number of active degrees of freedom, while in Figure \ref{fig:unif_2_parts} we report the various components of the a posteriori estimator, considering two cases with $Re=1$ and $Re=40$. We see that in this test $\eta_f$ is different from $0$ but decreases fast. Furthermore, we notice that for $Re=1$ the dominant part is again the bulk error, whereas for $Re=40$ at the beginning the dominant part is $\eta_{c2}$, then again $\eta_B$ becomes the most important one. Again, Newton's method converges in 3 iterations or less. Overall, velocity and pressure errors as well as the estimator scale according to the optimal convergence rate.

\subsection{Fluid dynamics case}
The fluid dynamics case describes the flow around a square cylinder in a channel, following the setting in \cite{Fluid_1}. There is no external forces, so $\mathbf{f}=\mathbf{0}$. At the inflow we consider a parabolic velocity profile in the $x$-component:
\begin{equation}
    \mathbf{u}_I=\begin{bmatrix}
-1/16(y-4)(y+4)\\
0
\end{bmatrix}\,.
\end{equation}
No slip conditions ($\mathbf{u}=\mathbf{0}$) are enforced on the cylinder and channel walls, whereas  homogeneous Neumann conditions are chosen at the outflow. 
The problem is adimensionalized with the maximum inflow velocity as characteristic velocity.

\begin{figure}
        \centering
        \begin{subfigure}[b]{0.5\textwidth}
                \centering
 \includegraphics[width=\textwidth]{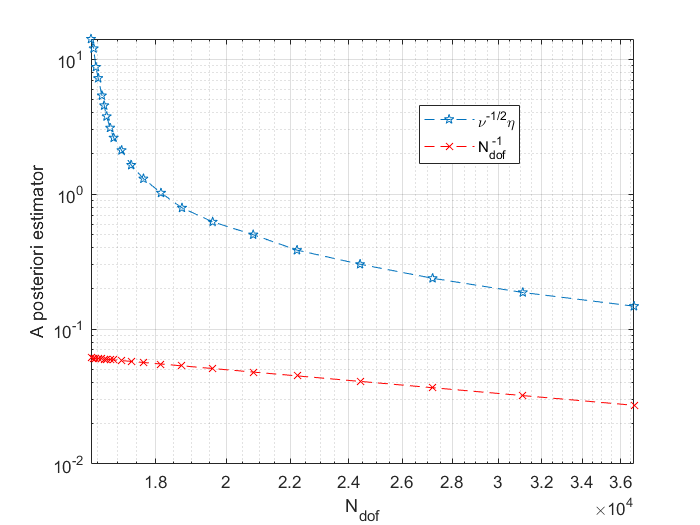}
         \caption{A posteriori estimator}
                \label{fig:sensor3}
        \end{subfigure}%    <-- % added here
        \hfill %% useful if width of each figure is less the .5\textwidth
        \begin{subfigure}[b]{0.5\textwidth}
                \centering
         \includegraphics[width=\textwidth]{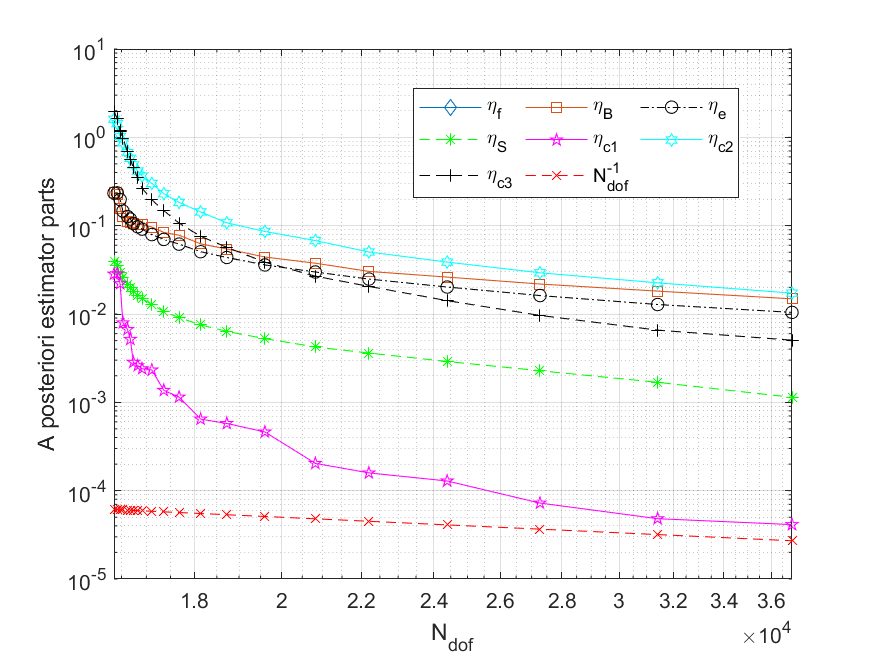}
         \caption{Components  of the estimator}
         \label{}
        \end{subfigure}
        \caption{Fluid dynamics case. A posteriori estimator (left) and its components  (right) for $Re=30$} \label{fig:estim_30}
\end{figure}

\begin{figure}[t!]
     \centering
     \begin{adjustbox}{minipage=0.9\linewidth,scale=1.1}
     \begin{subfigure}[b]{0.49\textwidth}
         \centering
         \includegraphics[width=\textwidth]{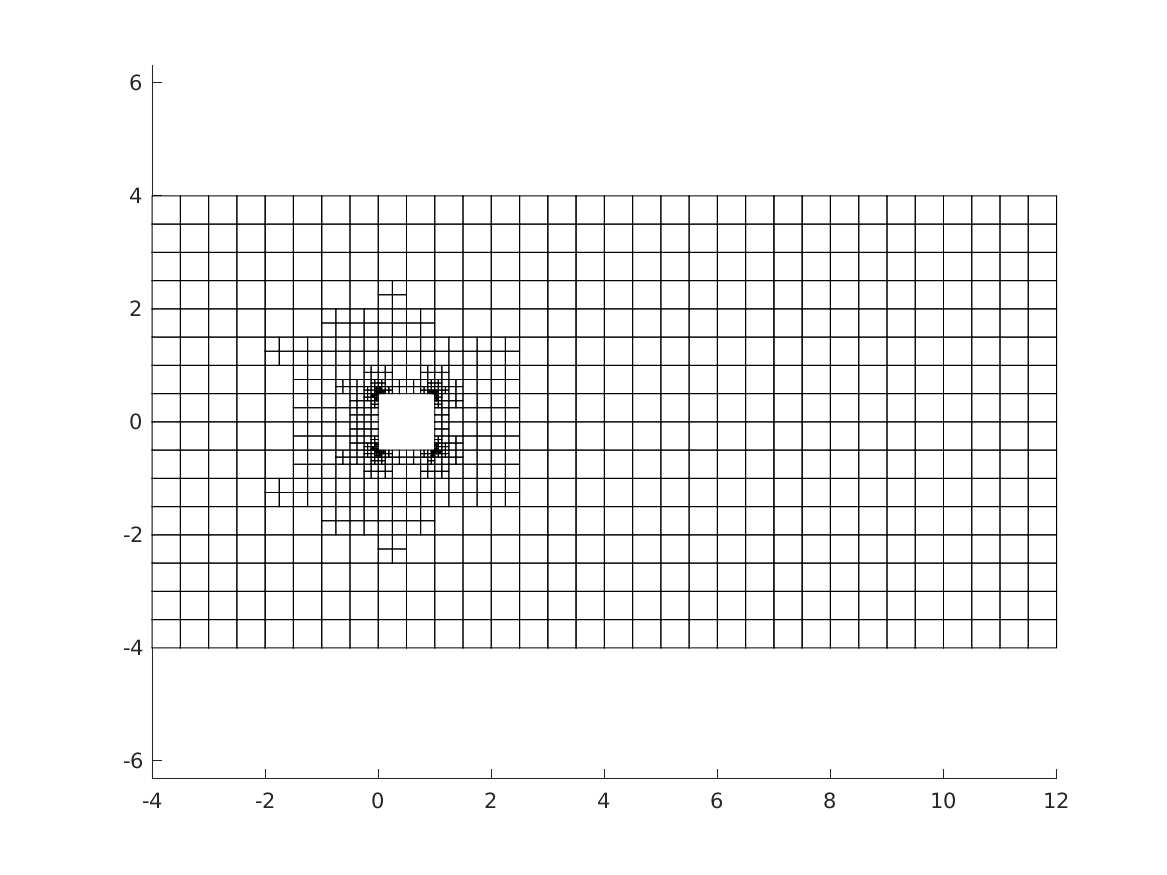}
         \caption{$Re$=1 (10 steps)}
         \label{}
     \end{subfigure}
     \hfill
    %  \begin{subfigure}[b]{0.45\textwidth}
    %      \centering
    %      \includegraphics[width=\textwidth]{}
    %      \caption{$Re$=5}
    %      \label{}
    %  \end{subfigure}
     \hfill
     \begin{subfigure}[b]{0.49\textwidth}
         \centering
         \includegraphics[width=\textwidth]{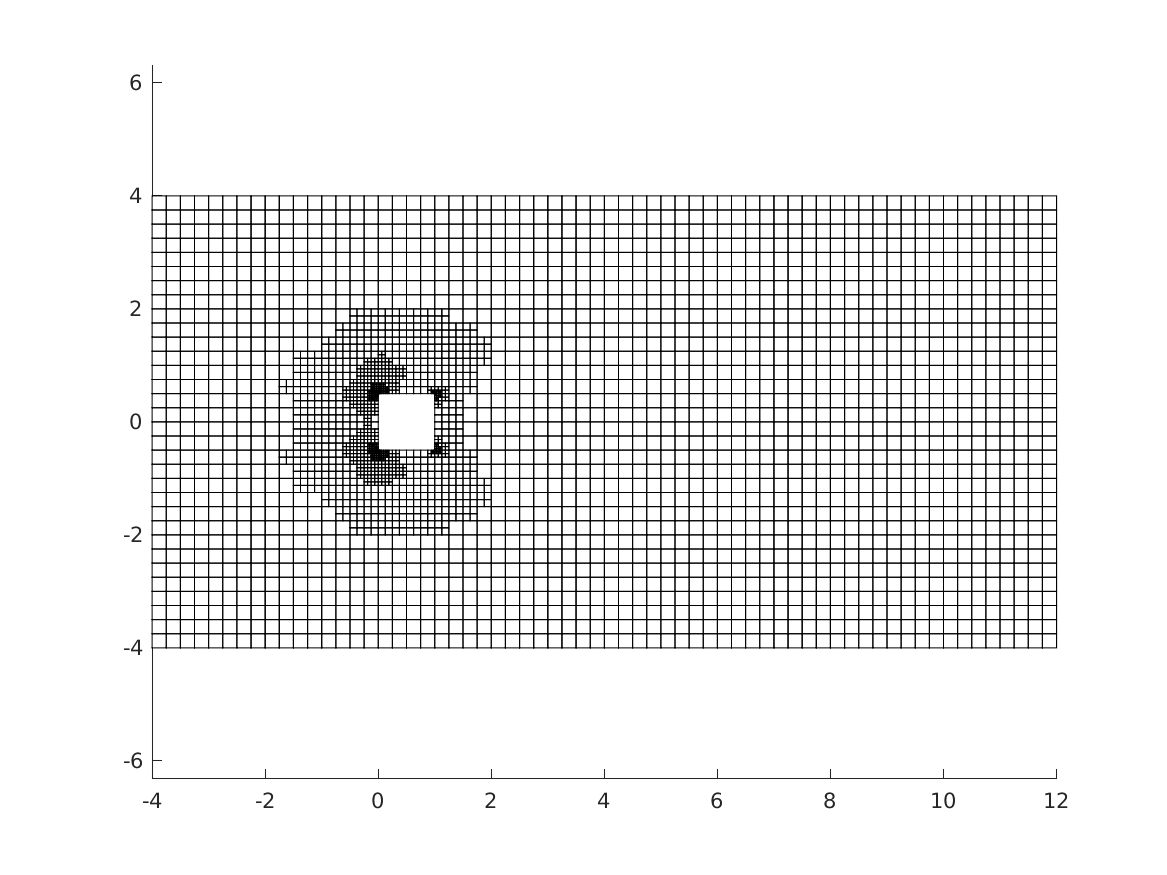}
         \caption{$Re$=10 (15 steps)}
         \label{}
     \end{subfigure}
     \hfill
     \begin{subfigure}[b]{0.49\textwidth}
         \centering
         \includegraphics[width=\textwidth]{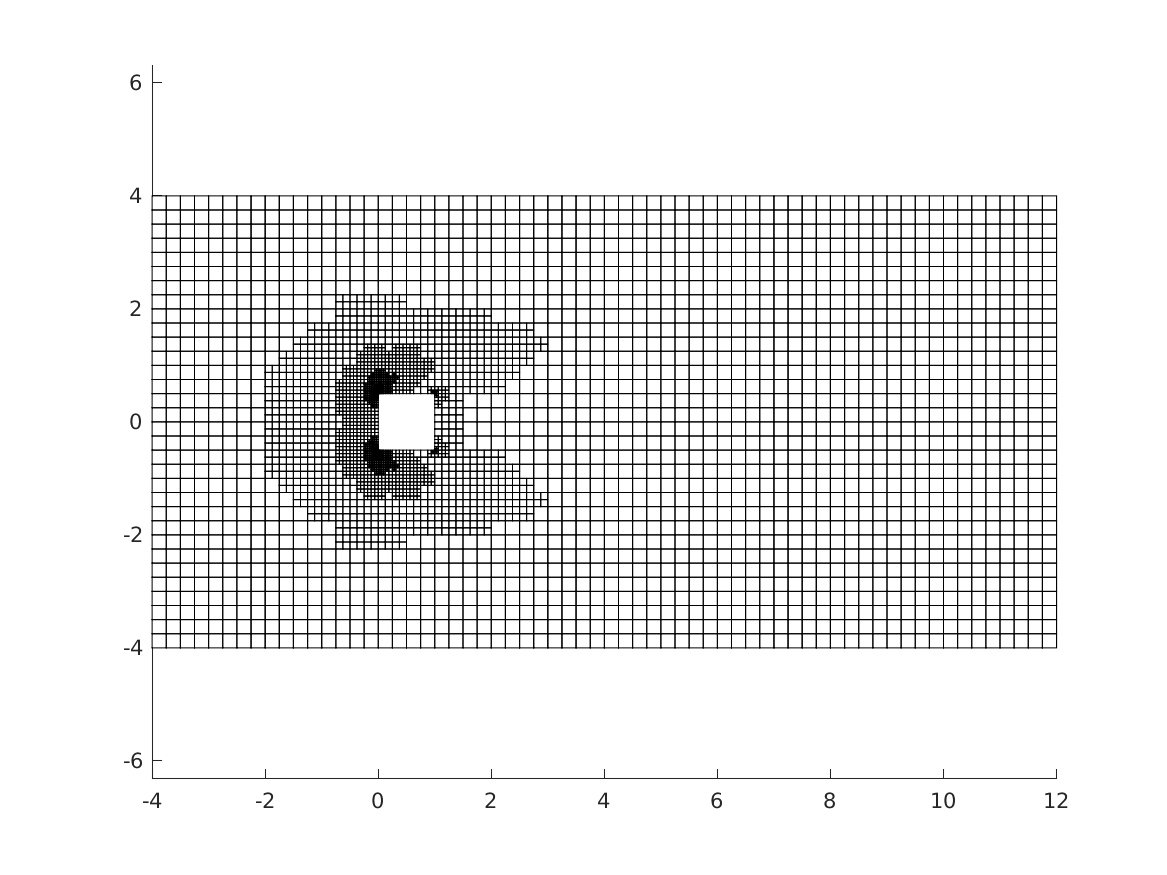}
         \caption{$Re$=30 (20 steps)}
         \label{fig:mesh_30}
     \end{subfigure}
    %   \begin{subfigure}[b]{0.45\textwidth}
    %      \centering
    %      \includegraphics[width=\textwidth]{fluid/mesh_caso_2040_16.png}
    %      \caption{$Re$=40 (15 steps)}
    %      \label{}
    %  \end{subfigure}
      \begin{subfigure}[b]{0.49\textwidth}
         \centering
         \includegraphics[width=\textwidth]{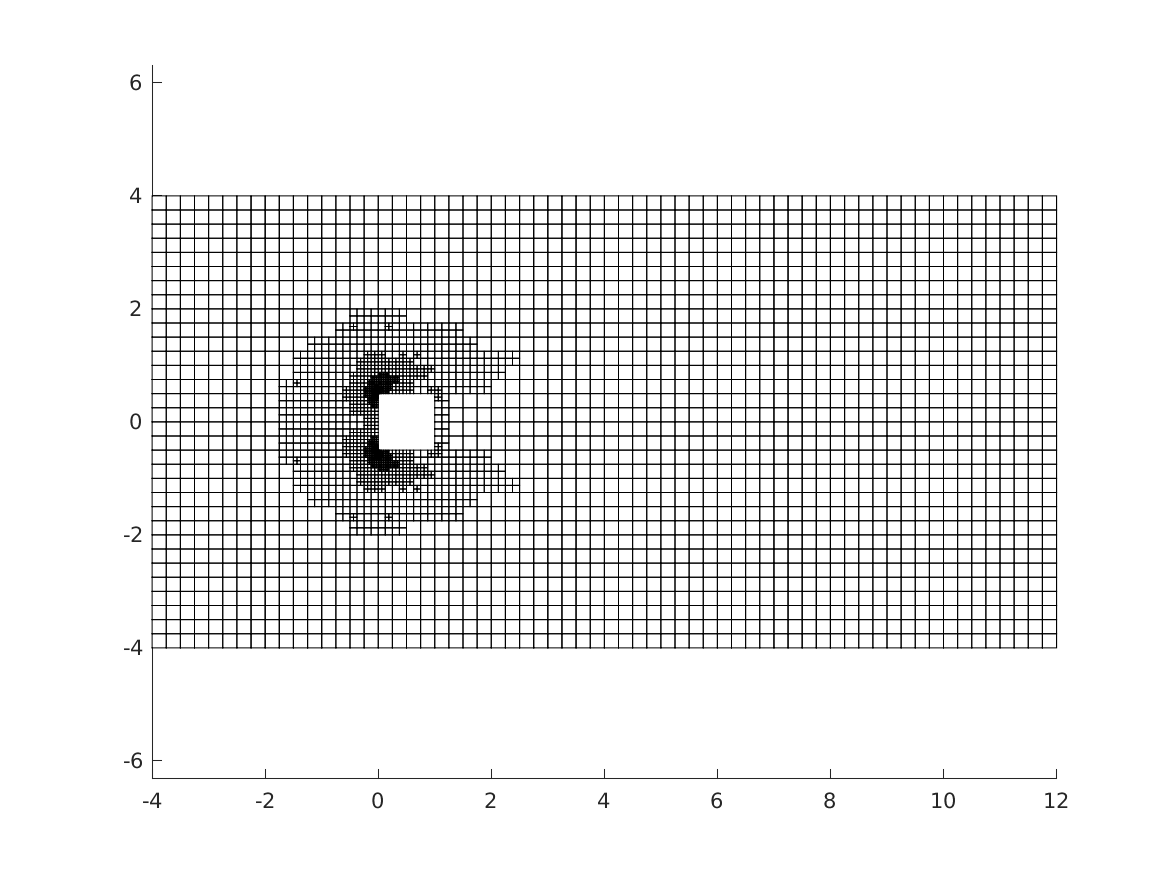}
         \caption{$Re$=50 (20 steps)}
         \label{}
     \end{subfigure}
     \end{adjustbox}
        \caption{Flow around a square cylinder. Meshes at the end of the adaptive refinement for different Reynolds numbers. `Steps' indicates the number of steps performed in the refinement loop}
        \label{fig:mesh_fluid}
\end{figure}

\begin{figure}[t!]
     \centering
     \begin{adjustbox}{minipage=0.9\linewidth,scale=1.1}
     \begin{subfigure}[b]{0.49\textwidth}
         \centering
         \includegraphics[width=\textwidth]{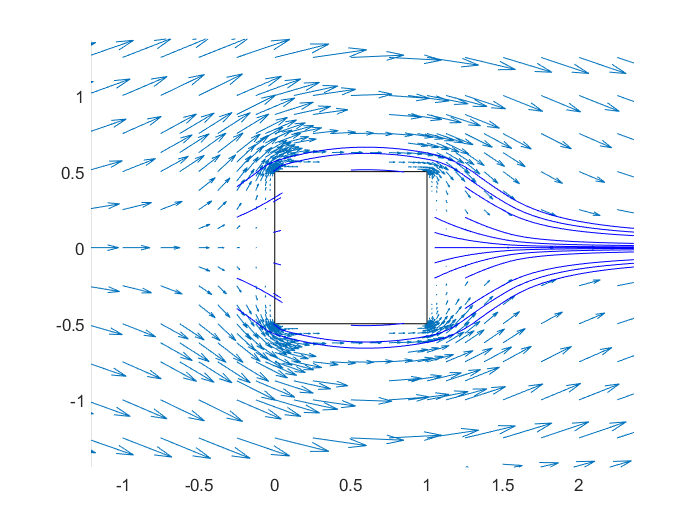}
         \caption{$Re$=1}
         \label{}
     \end{subfigure}
    %  \hfill
    %  \begin{subfigure}[b]{0.45\textwidth}
    %      \centering
    %      \includegraphics[width=\textwidth]{fluid/errori_caso_7001.png}
    %      \caption{$Re$=5}
    %      \label{}
    %  \end{subfigure}
     \hfill
     \begin{subfigure}[b]{0.49\textwidth}
         \centering
         \includegraphics[width=\textwidth]{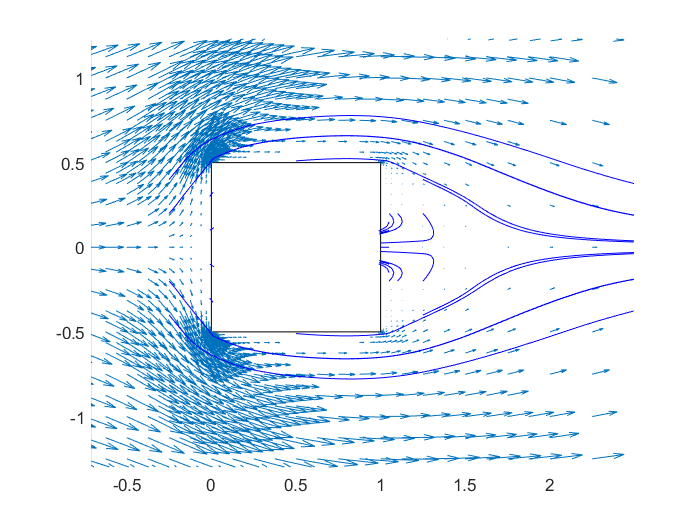}
         \caption{$Re$=10}
         \label{}
     \end{subfigure}
     \hfill
     \begin{subfigure}[b]{0.49\textwidth}
         \centering
         \includegraphics[width=\textwidth]{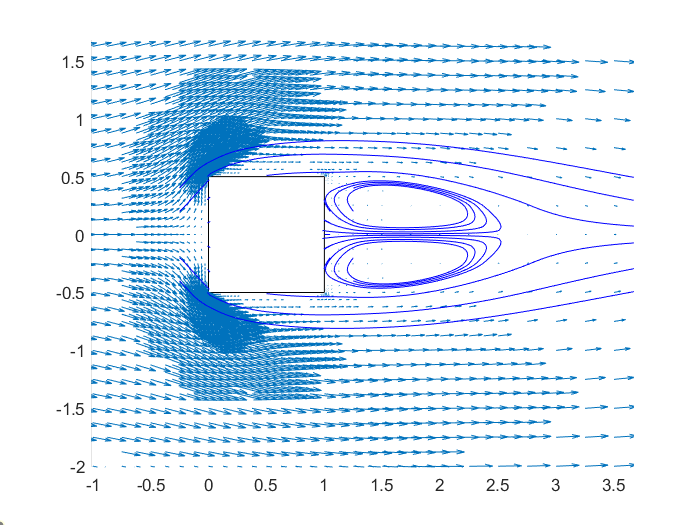}
         \caption{$Re$=30}
         \label{fig:quiver_30}
     \end{subfigure}
    %   \begin{subfigure}[b]{0.45\textwidth}
    %      \centering
    %      \includegraphics[width=\textwidth]{fluid/errori_caso_7004.png}
    %      \caption{$Re$=40}
    %      \label{}
    %  \end{subfigure}
      \begin{subfigure}[b]{0.49\textwidth}
         \centering
         \includegraphics[width=\textwidth]{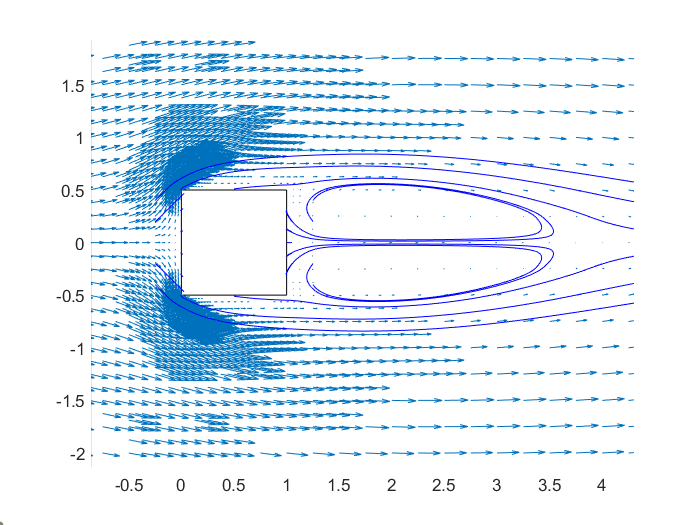}
         \caption{$Re$=50}
         \label{}
     \end{subfigure}
      \end{adjustbox}
        \caption{Flow around a square cylinder. Velocity fields and streamlines, computed on the corresponding meshes given in Figure \ref{fig:mesh_fluid}}
        \label{fig:quiver_fluid}
\end{figure}

 We apply in this case an adaptive refinement with D\"orfler parameter $\theta=0.4$. 
In Figure \ref{fig:estim_30} we report the estimator and its components vs the number of active degrees of freedom, when $Re=30$: we observe that all monitored quantities decrease very fast at the beginning, while later they follow the expected convergence. We also show in Figure \ref{fig:mesh_fluid} some meshes produced by the refinement, and in Figure \ref{fig:quiver_fluid} some velocity fields, where the blue lines represent \textit{Matlab}-generated approximate streamlines. One may appreciate the refinement around the cylinder, in particular around its corners pointing towards the inflow side, as we can see from Figure \ref{fig:mesh_30}.  Indeed,  the refinement  becomes, as the Reynolds number increases, more and more non-symmetric with respect to the vertical axis $x=1/2$ passing through the center of the cylinder (whereas in the Stokes case the refined meshes are symmetric in the flow direction). The refinement is concentrated about the upstream edge of the cylinder as $Re$ increases, while there is less refinement near the downstream edge.

In the current situation, we observe the appearance of recirculation eddies past the cylinder (see Figure \ref{fig:quiver_30}), whose length (the {\em recirculation length}) can be estimated by linear interpolation, after identifying the element edge on the line $y=0$ where $u_x$ (the \textit{x}-component of $\mathbf{u}$) changes sign past the cylinder.  The computed values of the recirculation length for different Reynolds numbers are reported in Table \ref{tab:table_re}; they exhibit a linear growth, which is in good agreement with analogous results in the literature (see \cite{Fluid_1}). Note that in the recirculation area we do not see a strong refinement; this can be related to the presence of very low values of velocity therein. 

% ????However, investigating the problem for different $Re<60$ (stationary and laminar range) we see that the boundary layer separates at a value of $Re$ in agreement with \cite{Fluid_1}, while the recicurlation length increases linearly with respect with the Reynolds number as appearing in the simulations in \cite{Fluid_1}.???? (INSERIRE?)

\begin{table}[h]
\caption{Recirculation length for different Reynolds numbers} %title of the table
\centering % centering table
\begin{tabular}{c c c c c c c c c c c} % creating eight columns

\hline % inserts single-line
\textit{Re} & 10 & 15 & 20& 25& 30& 35& 40 &45 &50 &55\\ % Entering row contents
Recir. Len. & 1.50 & 1.79  & 2.09& 2.39&2.70 & 3.00& 3.29 & 3.58&3.79 &4.15\\

\hline % inserts single-line
\end{tabular}
\label{tab:table_re}
\end{table}


\section{Conclusions}
\label{sec:conclusions}
In this work we considered the divergence-free Virtual Element Method proposed in \cite{Navier-Stokes_1} for the numerical solution of the steady incompressible Navier-Stokes equations. The method is of arbitrary order $k \geq 2$, meaning that the local VEM space approximates velocities on the boundary of the mesh element by piecewise polynomials of degree $k$, whereas pressures are locally approximated by polynomials of degree $k-1$. 
\begin{itemize}
    \item We defined a residual-based a posteriori estimator for a linear combination of the $H^1$-seminorm error on velocity and the $L^2$-norm error on pressure. 
    \item We established the reliability of the new estimator, proving that it provides an upper bound for this error. Furthermore, we highlighted the origin of the different components of the estimator. 
    \item With respect to the Stokes case, in addition to a different definition of the bulk error, our estimator includes three new terms originated by the VEM discretization of the nonlinear convective term.
    \item We derived lower bounds for the error, which involve all the terms in the estimator and indicate that the proposed estimator is also efficient. 
    \item We performed some numerical tests with order $k=2$ and uniform mesh refinement, showing the correct rate of decay of the estimator and its control over the true errors. We also highlighted the behaviour of the various components of the estimator.
    \item Finally, we applied an adaptive mesh refinement based on the new estimator to the flow around a square obstacle in the stationary and laminar regime $Re<60$. Our results are in good agreement with classical results obtained by finite-volume simulations (see \cite{Fluid_1}); in particular, the recirculation lengths past the obstacle increase in agreement with \cite{Fluid_1}. The refined meshes appear to fit the structures of the physical solution for simulations at different Reynolds numbers, showing a good performance of the new estimator.
\end{itemize}

\bibliographystyle{unsrtnat}
%\bibliography{references}  %%% Uncomment this line and comment out the ``thebibliography'' section below to use the external .bib file (using bibtex) .





% \end{thebibliography}

\section*{Declarations}
The authors have no relevant financial or non-financial interests to disclose.

This research was done in the framework of the Italian MIUR Award ``Dipartimenti di Eccellenza 2018-2022" granted to the Department of Mathematical Sciences, Politecnico di Torino (CUP: E11G18000350001). CC is a member of the Italian INdAM-GNCS research group and was supported by the MIUR PRIN Project 201752HKH8-003.

\bibliography{Canuto_Rosso_apost_VEM-bibliography.bib}% common bib file


\end{document}